\documentclass[10pt]{amsart}

\numberwithin{equation}{section} 

\newcommand{\bC}{\mathbb{C}}
\newcommand{\bZ}{\mathbb{Z}}

\newcommand{\ch}{\mathrm{ch}}

\newcommand{\ord}{{\mathrm{ord}}}
\newcommand{\Sol}{{\mathrm{Sol}}}
\newcommand{\Res}{\operatornamewithlimits{Res}}
\newcommand{\cJ}{\mathcal{J}}
\newcommand{\cN}{\mathcal{N}}
\newcommand{\cH}{\mathcal{H}}
\newcommand{\cL}{\mathcal{L}}

\newcommand{\cR}{\mathcal{R}}
\newcommand{\cK}{\mathcal{K}}
\newcommand{\cW}{\mathcal{W}}
\newcommand{\Uq}{U_q(X_n^{(1)})}
\newcommand{\Uqcl}{U_q(X_n)}
\newcommand{\tQ}{\tilde{Q}}

\newcommand{\tR}{\tilde{R}}
\newcommand{\tK}{\tilde{K}}
\newcommand{\tx}{\tilde{x}}

\newcommand{\sa}{^{(a)}}
\newcommand{\sab}{^{(b)}}
\newcommand{\am}{_m^{(a)}}
\newcommand{\bk}{_k^{(b)}}
\newcommand{\XXX}{$X\hskip-1.2pt X\hskip-1.2pt X$}
\newcommand{\XXZ}{$X\hskip-1.2pt X\hskip-1.2pt Z$}

\usepackage{amsthm}
\newtheorem{thm}{Theorem}[section]
\newtheorem{lem}[thm]{Lemma}
\newtheorem{prop}[thm]{Proposition}
\newtheorem{cor}[thm]{Corollary}
\newtheorem{conj}[thm]{Conjecture}

\newtheorem{thmdef}[thm]{Theorem-Definition}
\theoremstyle{definition}
\newtheorem{defn}[thm]{Definition}
\newtheorem{exmp}[thm]{Example}
\theoremstyle{remark}
\newtheorem{rem}[thm]{Remark}

\begin{document}

\title[Bethe Equation and weight multiplicities]
{Bethe Equation at $q=0$, M\"obius
 Inversion
Formula, and Weight Multiplicities:\\ II. $X_n$ case}

\author{Atsuo Kuniba}
\address{
Institute of Physics, University of Tokyo,
Tokyo 153-8902, Japan}
\curraddr{}
\email{atsuo@gokutan.c.u-tokyo.ac.jp}
%
\author{Tomoki Nakanishi}
\address{Graduate School of Mathematics,
    Nagoya University, Nagoya 464-8602, Japan}
\email{nakanisi@math.nagoya-u.ac.jp}
%
%

\begin{abstract}
We study a family of power series  characterized by
a system of recursion relations
 ($Q$-system) with a certain convergence property.
We show that the coefficients of the series
are expressed by the numbers which formally
count
the off-diagonal solutions of
 the $U_q(X^{(1)}_n)$ Bethe equation at $q=0$.
The series are conjectured to be the
$X_n$-characters of a certain family of
irreducible finite-dimensional $U_q(X^{(1)}_n)$-modules
which we call the KR (Kirillov-Reshetikhin) modules.
Under the above conjecture, 
these coefficients give a formula of
 the weight multiplicities of the tensor products of
the KR modules, which is also interpreted as the
formal completeness of the \XXZ-type Bethe vectors.
\end{abstract}

\maketitle


\section{Introduction}
\label{sec:intro}

\subsection{Background}

\par Since the early days of the
algebraic Bethe ansatz approach to the 
integrable lattice models,
the finite-dimensional representations of 
the affine quantum group
$U_q(X_n^{(1)})$ and its sister,
Yangian $Y(X_n)$, have attracted 
much attention.
The categories of finite-dimensional modules of
$U_q(X_n^{(1)})$ and  $Y(X_n)$ are equivalent \cite{CP,D}.
The case $X_n=A_n$ is rather well understood \cite{CP,Ar}.
Also, a general theory of 
the character has been constructed \cite{FR,Kn}.
However, no universal description of the character,
 like the Weyl character formula,
is known so far.
\par
In the series of works
\cite{K1,K2,K3,KR},
Kirillov and Reshetikhin focused
 attention on a special class
of $Y(X_n)$-modules.
We call them (and the corresponding $U_q(X_n^{(1)})$-modules)
the {\em KR modules}.
The KR modules of 
$Y(X_n)$ (resp.\ $U_q(X_n^{(1)})$)
provide natural generalizations of
the well-known \XXX\ (resp.\ \XXZ) spin chain.
They found that the formal counting of
the Bethe vectors of the \XXX-type spin chains,
with the hypothesis of the completeness of the
Bethe vectors,
leads to a remarkable conjectural formula for
 the {\em multiplicity of the $X_n$-irreducible 
components\/}
in  the tensor products of the 
KR modules of $Y(X_n)$.

\par In this paper
we show an analogous phenomenon occurs also for the \XXZ-type case.
In contrast with the \XXX-type case, however,
the formal counting
of the Bethe vectors of the \XXZ-type spin chains
 now leads to a conjectural formula
for the {\em weight multiplicity\/} in
the tensor products of the KR modules of $\Uq$.
We hope that these formulae
will guide us to the proper understanding of 
the KR modules.
This paper is the continuation of Part I \cite{KN}, where the case
$X_n=A_1$ is described.
Below we shall formulate 
and explain our problem more precisely.

\subsection{KR modules}
\label{subsec:rect}

\par Let $X_n$ be one of the finite-dimensional
 simple Lie algebras over $\mathbb{C}$
and $U_q(X_n^{(1)})$ be the non-twisted
quantum affine algebra
 associated to
$X_n$ without the derivation operator.
Let $\alpha_a$ and $\Lambda_a$ ($a=1,\dots,n$) be the
 simple roots and
fundamental weights of $X_n$.
We enumerate the vertices of the Dynkin diagram
as in {Table} \ref{tab:Dynkin}.
Let $(\,\cdot\,|\,\cdot\,)$ be the standard bilinear
form normalized as $(\alpha|\alpha)=2$ for a long root $\alpha$.
We set
$t_a   =   {2}/{(\alpha_a \vert \alpha_a)} \in \{\, 1, 2, 3\, \}$.
The Cartan matrix is $
C_{a b}  =  t_a(\alpha_a \vert \alpha_b)$,
and
$\alpha_a  =  \sum_{b=1}^n C_{b a} \Lambda_b$.

\unitlength=0.8pt
\begin{table}
\caption{Dynkin diagrams}
\label{tab:Dynkin}
\begin{tabular}[t]{rlrlrl}
$A_n$:&
\begin{picture}(106,20)(-5,-5)
\multiput( 0,0)(20,0){2}{\circle{6}}
\multiput(80,0)(20,0){2}{\circle{6}}
\multiput( 3,0)(20,0){2}{\line(1,0){14}}
\multiput(63,0)(20,0){2}{\line(1,0){14}}
\multiput(39,0)(4,0){6}{\line(1,0){2}}
\put(0,-8){\makebox(0,0)[t]{$1$}}
\put(20,-8){\makebox(0,0)[t]{$2$}}
\put(80,-8){\makebox(0,0)[t]{$n\!\! - \!\! 1$}}
\put(100,-11){\makebox(0,0)[t]{$n $}}
\end{picture}
&\quad $E_6$:&
\begin{picture}(86,40)(-5,-5)
\multiput(0,0)(20,0){5}{\circle{6}}
\put(40,20){\circle{6}}
\multiput(3,0)(20,0){4}{\line(1,0){14}}
\put(40, 3){\line(0,1){14}}
\put( 0,-8){\makebox(0,0)[t]{$1$}}
\put(20,-8){\makebox(0,0)[t]{$2$}}
\put(40,-8){\makebox(0,0)[t]{$3$}}
\put(60,-8){\makebox(0,0)[t]{$4$}}
\put(80,-8){\makebox(0,0)[t]{$5$}}
\put(48,20){\makebox(0,0)[l]{$6$}}
\end{picture}
& $F_4$:&
\begin{picture}(66,20)(-5,-5)
\multiput( 0,0)(20,0){4}{\circle{6}}
\multiput( 3,0)(40,0){2}{\line(1,0){14}}
\multiput(22.85,-1)(0,2){2}{\line(1,0){14.3}}
\put(30,0){\makebox(0,0){$>$}}
\put(0,-8){\makebox(0,0)[t]{$1$}}
\put(20,-8){\makebox(0,0)[t]{$2$}}
\put(40,-8){\makebox(0,0)[t]{$3$}}
\put(60,-8){\makebox(0,0)[t]{$4$}}
\end{picture}
\\
$B_n$:&
\begin{picture}(106,20)(-5,-5)
\multiput( 0,0)(20,0){2}{\circle{6}}
\multiput(80,0)(20,0){2}{\circle{6}}
\multiput( 3,0)(20,0){2}{\line(1,0){14}}
\multiput(63,0)(20,0){1}{\line(1,0){14}}
\multiput(82.85,-1)(0,2){2}{\line(1,0){14.3}}
\multiput(39,0)(4,0){6}{\line(1,0){2}}
\put(90,0){\makebox(0,0){$>$}}
\put(0,-8){\makebox(0,0)[t]{$1$}}
\put(20,-8){\makebox(0,0)[t]{$2$}}
\put(80,-8){\makebox(0,0)[t]{$n\!\! -\!\! 1$}}
\put(100,-11){\makebox(0,0)[t]{$n$}}
\end{picture}
& $E_7$:&
\begin{picture}(106,40)(-5,-5)
\multiput(0,0)(20,0){6}{\circle{6}}
\put(40,20){\circle{6}}
\multiput(3,0)(20,0){5}{\line(1,0){14}}
\put(40, 3){\line(0,1){14}}
\put( 0,-8){\makebox(0,0)[t]{$1$}}
\put(20,-8){\makebox(0,0)[t]{$2$}}
\put(40,-8){\makebox(0,0)[t]{$3$}}
\put(60,-8){\makebox(0,0)[t]{$4$}}
\put(80,-8){\makebox(0,0)[t]{$5$}}
\put(100,-8){\makebox(0,0)[t]{$6$}}
\put(48,20){\makebox(0,0)[l]{$7$}}
\end{picture}
& $G_2$:&
\begin{picture}(26,20)(-5,-5)
\multiput( 0, 0)(20,0){2}{\circle{6}}
\multiput(2.68,-1.5)(0,3){2}{\line(1,0){14.68}}
\put( 3, 0){\line(1,0){14}}
\put( 0,-8){\makebox(0,0)[t]{$1$}}
\put(20,-8){\makebox(0,0)[t]{$2$}}
\put(10, 0){\makebox(0,0){$>$}}
\end{picture}
\\
$C_n$:&
\begin{picture}(106,20)(-5,-5)
\multiput( 0,0)(20,0){2}{\circle{6}}
\multiput(80,0)(20,0){2}{\circle{6}}
\multiput( 3,0)(20,0){2}{\line(1,0){14}}
\multiput(63,0)(20,0){1}{\line(1,0){14}}
\multiput(82.85,-1)(0,2){2}{\line(1,0){14.3}}
\multiput(39,0)(4,0){6}{\line(1,0){2}}
\put(90,0){\makebox(0,0){$<$}}
\put(0,-8){\makebox(0,0)[t]{$1$}}
\put(20,-8){\makebox(0,0)[t]{$2$}}
\put(80,-8){\makebox(0,0)[t]{$n\!\! -\!\! 1$}}
\put(100,-11){\makebox(0,0)[t]{$n$}}
\end{picture}
&
$E_8$:&
\begin{picture}(126,40)(-5,-5)
\multiput(0,0)(20,0){7}{\circle{6}}
\put(80,20){\circle{6}}
\multiput(3,0)(20,0){6}{\line(1,0){14}}
\put(80, 3){\line(0,1){14}}
\put( 0,-8){\makebox(0,0)[t]{$1$}}
\put(20,-8){\makebox(0,0)[t]{$2$}}
\put(40,-8){\makebox(0,0)[t]{$3$}}
\put(60,-8){\makebox(0,0)[t]{$4$}}
\put(80,-8){\makebox(0,0)[t]{$5$}}
\put(100,-8){\makebox(0,0)[t]{$6$}}
\put(120,-8){\makebox(0,0)[t]{$7$}}
\put(88,20){\makebox(0,0)[l]{$8$}}
\end{picture}
\\
$D_n$:&
\begin{picture}(106,40)(-5,-5)
\multiput( 0,0)(20,0){2}{\circle{6}}
\multiput(80,0)(20,0){2}{\circle{6}}
\put(80,20){\circle{6}}
\multiput( 3,0)(20,0){2}{\line(1,0){14}}
\multiput(63,0)(20,0){2}{\line(1,0){14}}
\multiput(39,0)(4,0){6}{\line(1,0){2}}
\put(80,3){\line(0,1){14}}
\put(0,-8){\makebox(0,0)[t]{$1$}}
\put(20,-8){\makebox(0,0)[t]{$2$}}
\put(80,-8){\makebox(0,0)[t]{$n\!\! - \!\! 2$}}
\put(103,-8){\makebox(0,0)[t]{$n\!\! -\!\! 1$}}
\put(85,20){\makebox(0,0)[l]{$n$}}
\end{picture}
\\
\end{tabular}
%
\end{table}


\par 
The irreducible finite-dimensional $\Uq$-modules are
parameterized by $n$-tuples of polynomials
({\em Drinfeld polynomials\/})   $(P_b)_{b=1}^n$
 with constant terms 1 \cite{CP}.
Here we follow the convention in \cite{FR}.

\begin{defn} 
For each $a\in \{\, 1,\dots,n\, \}$, $m \in \{\, 1,2,\dots\,\}$,
and $u\in \mathbb{C}$,
the irreducible finite-dimensional 
$\Uq$-module whose  Drinfeld
polynomials $(P_b)$ are
\begin{align}\label{eq:rect1}
P_b(v)=
\begin{cases}
 \prod_{j=1}^m (1-vq^u q^{(m-2j+1)/t_a})
 &b=a\\
1& \text{otherwise}
\end{cases}
\end{align}
is called a {\em KR (Kirillov-Reshetikhin) module\/} and denoted by
$W\am(u)$.
\end{defn}

\par 
For $m=1$, the KR modules are the fundamental modules.
Through the injection $\Uqcl \rightarrow
\Uq$, any $\Uq$-module $W$ is  regarded as
a $\Uqcl$-module.

\begin{lem}\label{lem:equiv1}
For any $u,\ u'\in \mathbb{C}$,
$W\am(u)$ and $W\am(u')$ are equivalent as
 $\Uqcl$-modules.
\end{lem}
\begin{proof}
It is well known that $W\am(u')$ is
obtained from $W\am(u)$ by a pull-back of
an automorphism $\sigma$ of $\Uq$ which
preserves $\Uqcl$.
\end{proof}

\par
In view of Lemma \ref{lem:equiv1},
let us write  the 
common underlying $\Uqcl$-module for
the family $W\am(u)$ as $W\am$.
In general, there is a natural identification of a $\Uqcl$-module $V$
with an $X_n$-module
through the limit $q\rightarrow 1$.
With this identification,
we use the $X_n$-weight and $X_n$-character to describe $V$,
instead of $U_q(X_n)$-weight and $U_q(X_n)$-character.
For example, $\ch\, W\am$ means the $X_n$-character of $W\am$.
The $X_n$-weight of the
 highest weight vector of $W\am(u)$ is $m\Lambda_a$.

\subsection{The series $Q\am$ and Kirillov-Reshetikhin's conjecture}

\par 
Let 
 $x_a=e^{\Lambda_a}$ and
$y_a=e^{-\alpha_a}$, $a=1$, \dots, $n$,
be the formal exponentials of
the fundamental weights and (the minus of)
the simple roots of $X_n$.
With the multivariable notation
$x=(x_a)_{a=1}^n$ and $y=(y_a)_{a=1}^n$,
we write the relation
$y_a=\prod_{b=1}^n x_b^{-(\alpha_a|\alpha_b)t_b}$
as $y=y(x)$.
Its inverse $x=x(y)$
 involves fractional powers for some $X_n$.
Let $\mathbb{C}[[y]]$
be the ring of (formal)
 power series of $y=(y_a)_{a=1}^n$ with the standard topology.
For a power series $f(y)$,
$f(y(x))$ is a Laurent series of $x$.
Conversely, for a Laurent series
$f(x)$, $f(x(y))$ is a 
Puiseax (fractional Laurent)
series  of $y$, in general.
\par
Let
$H$ denote the index set
\begin{align}\label{eq:hdef1}
H=\{\, (a,m)\mid
a\in \{1,\dots,n\},\ m \in \{1,2,\dots\} \,\},
\end{align}
and for $(a,m)$, $(b,k)\in H$, we define
\begin{equation}\label{eq:bfunc}
B_{am,bk}=
2\min(t_bm,t_ak)-\min(t_bm,t_a(k+1))
-\min(t_bm,t_a(k-1)).
\end{equation}
\par

The following series are our main interest in this paper.
\begin{thmdef}\label{defn:unique}
Let $(\tQ\am(y))_{(a,m)\in H}$ be
 the unique family  of the invertible
 power series of $y$
which satisfies  (\~Q-I) and (\~Q-II):
\par (\~{Q}-I). ($Q$-system) Let $\tQ_0\sa(y) =1$.
For $m=1$, $2$, \dots,
\begin{equation}
\label{eq:qsys6}
(\tQ\am(y))^2=
\tQ_{m+1}\sa(y) \tQ_{m-1}\sa(y)
+
y_a^m
(\tQ\am(y))^2
\prod_{\scriptstyle
(b,k)\in H}
(\tQ\bk(y))^{-(\alpha_a|\alpha_b)B_{am,bk}}.
\end{equation}
\par (\~Q-II).
(convergence property)
The limit\/ $\lim_{m\to \infty} \tQ\am(y)$ exists
in $\mathbb{C}[[y]]$.
\par\noindent
Let $Q\am(x):=x_a^m \tQ\am(y(x))$.
Equivalently, 
 $(Q\am(x))_{(m,a)\in H}$ is the unique family of the Laurent
 series of $x$
which satisfies   (Q-I) and (Q-II):
\par (Q-I). ($Q$-system)
Let $Q_0\sa(x) =1$. For\/ $m=1$, $2$, \dots,

\begin{equation}
\label{eq:qsys1}
(Q\am(x))^2=
Q_{m+1}\sa(x) Q_{m-1}\sa(x)
+
(Q\am(x))^2
\prod_{\scriptstyle
(b,k)\in H}
(Q\bk(x))^{-(\alpha_a|\alpha_b)B_{am,bk}}.
\end{equation}
\par (Q-II).
(convergence property)
$\tQ\am(y):=x_a^{-m}
Q\am(x)|_{x=x(y)}$ is an invertible power series of $y$,
and the limit\/ $\lim_{m\to \infty} 
\tQ\am(y)$ exists in $\bC[[y]]$.
\par
\end{thmdef}

The existence and the uniqueness of $(\tQ\am(y))$,
together with two explicit expressions, will be
shown in Theorem \ref{thm:main3}, which is our
key theorem.

\begin{rem}\label{rem:Qsys}
The equivalence of the relations (\~Q-I) and (Q-I) is
easily seen with (\ref{eq:bfunc4}).
The product in the RHSs of
(\~Q-I) and (Q-I) are actually  finite products
(Proposition \ref{prop:bfunc1} (i)).
{}From given invertible power series
$\tQ_1^{(1)}$, \dots, $\tQ_1^{(n)}$,
the relation
(\~Q-I)  recursively determines all the
other invertible power series $\tQ\am$
(Proposition \ref{prop:rec1}).
\end{rem}

\begin{exmp}
For $A_1$, the relation (\ref{eq:qsys6}) becomes
\begin{equation}\label{eq:qsys8}
(\tQ^{(1)}_m(y))^2 =
\tQ^{(1)}_{m+1}(y)\tQ^{(1)}_{m-1}(y) + y_1^m.
\end{equation}
It is easy to check that $\tQ^{(1)}_{m}(y)=\sum_{j=0}^{m}
y_1^j$ satisfies the relation (\ref{eq:qsys8}) and that
$\lim_{m\to \infty}\tQ^{(1)}_{m}(y)=\sum_{j=0}^\infty y_1^j$.
Therefore, $Q^{(1)}_m(x)=x_1^m+x_1^{m-2}+\dots
+ x_1^{-m}$, which is the irreducible $A_1$-character
with highest weight $m\Lambda_1$.
\end{exmp}

Let $ \chi(\lambda)$ be
the irreducible $X_n$-character with highest weight
$\lambda$.
The following two theorems
were first proved  for $A_n$ by \cite{K2},
then generalized for the rest by \cite{HKOTY}:

\begin{thm}[\cite{K2,HKOTY}]\label{thm:hkoty1}
 For $X_n$ of classical type, i.e.,
$X_n=A_n$, $B_n$, $C_n$, and $D_n$,
$Q\am$'s in (\ref{eq:a1})--(\ref{eq:d1}) satisfy 
(Q-I) and (Q-II).
\begin{alignat}{2}
\label{eq:a1}
A_n:&&\quad
Q\am&= \chi(m\Lambda_a),\\
\label{eq:b1}
B_n:&&
Q\am &= \sum  \chi
(k_{a_0}\Lambda_{a_0}+k_{a_0+2}\Lambda_{a_0+2}+
\cdots +k_a\Lambda_a),\\
\label{eq:c1}
C_n:&&\quad
Q\am&=
\begin{cases}
\sum  \chi(k_{1}\Lambda_{1}
+k_{2}\Lambda_{2}+
\cdots +k_{a}\Lambda_{a})& 1\leq a \leq n-1\\
 \chi(m\Lambda_n)& a=n,
\end{cases}\\
\label{eq:d1}
D_n:&&
Q\am &=
\begin{cases}
 \sum  \chi
(k_{a_0}\Lambda_{a_0}+k_{a_0+2}\Lambda_{a_0+2}+
\cdots +k_a\Lambda_a)& 1\leq a \leq n-2\\
 \chi(m\Lambda_a)&a=n-1,n.\\
\end{cases}
\end{alignat}
In (\ref{eq:b1}) and (\ref{eq:d1}),
$a_0=0$ or $1$ with $a_0\equiv a$ {\rm mod} $2$, 
$\Lambda_0=0$,
and the sum is taken over non-negative integers
$k_{a_0}$, \dots, $k_a$ that satisfy
$t_a(k_{a_0}+\cdots +k_{a-2})+k_a=m$.
In (\ref{eq:c1}), the sum is taken over non-negative
integers $k_1$, \dots, $k_a$ that satisfy
$k_1+\cdots + k_a\leq m$, $k_b\equiv m \delta_{ab}$
{\rm mod} $2$.
\end{thm}
Because of the uniqueness, $Q\am(x)$ in Theorem \ref{thm:hkoty1}
coincides with the one in Definition
\ref{defn:unique}.
 In particular, $Q\am(x)$ is $\cW$ invariant for
any 
$X_n$ of classical type, where $\cW$ is the Weyl group of $X_n$.
\begin{thm}[\cite{K2,HKOTY}]\label{thm:hkoty2}
Let $X_n$ be arbitrary.
If $Q\am(x)$ is $\cW$ invariant for any $(a,m)\in H$,
then
\begin{align}\label{eq:qexp1}
Q\am&=\sum_{\lambda\in P_+}
\left\{\sum_{N\in \cN_\lambda} 
\prod_{(b,k)\in H}
\binom{P\bk+N\bk}{N\bk}
\right\} \chi(\lambda),
\end{align}
where $P_+$ is the set of the dominant integral
weights of $X_n$, and
\begin{align}
\label{eq:pkb1}
P\bk&=\min(k,m)\delta_{ab}
-\sum_{(c,j)\in H}(\alpha_b|\alpha_c)\min(t_ck,t_bj)N_j^{(c)},\\
\label{eq:nset5}
\cN_\lambda &=
\{\,N=(N\bk)_{(b,k)\in H}\mid 
N\bk \in \bZ_{\geq 0},\ 
m\Lambda_a-
\sum_{(b,k)\in H}k N\bk\alpha_b
=\lambda\,\}.
\end{align}
In particular, (\ref{eq:qexp1}) holds for 
any  $X_n$ of classical type.
\end{thm}

The following fundamental conjecture
is due to  Kirillov and Reshetikhin:
\begin{conj}[\cite{K3,KR}]\label{conj:KR1}
(i) For any   $X_n$ of classical type,
$\ch\, W\am$ equals to the RHSs of
(\ref{eq:a1})--(\ref{eq:d1}), respectively.
\par
(ii) For any $X_n$, $\ch\, W\am$ equals to
 the RHS of (\ref{eq:qexp1}).
\par
(iii) For any $X_n$, $\ch\, W\am$'s satisfy
the relation (Q-I) with $Q\am$ being replaced with
$\ch\, W\am$.

\end{conj}

\begin{rem}
Precisely speaking, 
 the existence of such modules $W\am(u)$
was claimed  in \cite{K3,KR} 
without the identification of their
Drinfeld polynomials (\ref{eq:rect1}).
According to Theorems \ref{thm:hkoty1} and \ref{thm:hkoty2},
for any  $X_n$ of classical type,
 (i) and (ii) are equivalent and (iii) follows from (i).
So far, Conjecture \ref{conj:KR1} has been completely
proved in the literature only for $A_n$ \cite{K2}
and $D_n$ \cite{Ch}.
\end{rem}
 
Since we expect that $Q\am$ are, in fact, $\cW$ invariant
also for any $X_n$ of exceptional type,
we reformulate Conjecture \ref{conj:KR1} 
more simply as
\begin{conj}\label{conj:qsys}
For any $X_n$, 
\begin{equation*}
 \ch\, W\am = Q\am.
\end{equation*}
\end{conj}

\subsection{Formal completeness of the \XXX-type Bethe vectors}
\label{subsec:problem}

Let us explain an interpretation of the expression (\ref{eq:qexp1})
in the context of  spin chains.
Let 
\begin{align}\label{eq:ndef1}
\cN=\{\, N=(N\am)_{(a,m)\in H}
\mid\, \text{$N\am\in \bZ_{\geq 0}$},\
\sum_{(a,m)\in H} N\am < \infty
\, \}.
\end{align}
For each $\nu=(\nu\am) \in \cN$,
 we associate a (finite) tensor product of
$W\am$'s,
\begin{equation}\label{eq:qspace}
\displaystyle W^\nu =
\bigotimes_{(a,m)\in H}
(W\am)^{\otimes \nu\am}.
\end{equation}
In the context of the spin chain, $W^\nu$ appears
as the {\em quantum space} on which the 
commuting family of the transfer matrices act.
Thus, we call $\nu$  the {\em quantum space data}.
For each $\nu\in \cN$,
we define a  Laurent series $Q^\nu(x)$ of $x$,
\begin{equation}\label{eq:Qnu}
Q^\nu(x)=\prod_{(a,m)\in H} (Q\am(x))^{\nu\am}
\end{equation}
and expand $Q^\nu(x)$ as
\begin{equation}\label{eq:Qexp4}
Q^\nu(x)
{\prod_{\alpha\in \Delta_+}(1-e^{-\alpha})}
=
\sum_{\lambda\in P} k^\nu_\lambda e^{\lambda},
\end{equation}
where $\Delta_+$ is the set of all the positive roots of $X_n$,
and $P$ is the weight lattice of $X_n$.
If Conjecture \ref{conj:qsys} is correct,
then  $Q^\nu=\ch\, W^\nu$.
It follows from the Weyl character formula that
$k^\nu_\lambda$ (for dominant $X_n$-weight $\lambda$)
is equal to the multiplicity 
$[W^\nu:V_\lambda]$ of the 
$\Uqcl$-irreducible components $V_\lambda$ with highest weight
$\lambda$
in $W^\nu$.
On the other hand, it was shown in \cite{K2,HKOTY}
 (cf.\ (\ref{eq:kwt2})) that 
\begin{equation}\label{eq:kres1}
k_\lambda^\nu = \sum_{N\in \cN_\lambda^\nu}
 K(\nu,N) 
\end{equation}
holds under the assumption of
the $\cW$ invariance for $Q\am$.
Here $K(\nu,N)$ and $\cN_\lambda^\nu$ are
 defined in (\ref{eq:Kdef})
and (\ref{eq:nset1}).
The RHS of (\ref{eq:kres1})
represents a formal counting 
of the Bethe vectors
of weight $\lambda$
for the  \XXX-type ($Y(X_n)$) spin chain with quantum space $W^\nu$
\cite{K1,K2,KR}.
We say the counting {\em formal\/}
because
$K(\nu,N)$ correctly counts the Bethe vectors 
only for special $(\nu,N)$'s.
Therefore, we call the equality (\ref{eq:kres1}) 
 the  {\em formal  completeness\/}
 of the \XXX-type 
Bethe vectors  (Corollary \ref{cor:completeness4}).
The formula (\ref{eq:qexp1}) is  the special case of
(\ref{eq:kres1}) for $\nu=(\nu\bk)$, $\nu\bk=\delta_{ab}\delta_{mk}$.

\subsection{Formal completeness of the \XXZ-type Bethe vectors}
\label{subsec:xxz}

\par
We remind that the most significant
difference between the Bethe vectors of the
spin chains of \XXX-type 
and \XXZ-type is that the former are $X_n$-singular,
while the latter are not necessarily $\Uqcl$-singular.
Accordingly, we expand $Q^\nu(x)$ as
\begin{equation}\label{eq:Qexp5}
Q^\nu(x)=
\sum_{\lambda\in P} r^\nu_\lambda e^\lambda.
\end{equation}
If Conjecture \ref{conj:qsys} is correct,
then $r^\nu_\lambda$ (for any $X_n$-weight $\lambda$)
is equal to
the weight multiplicity $\dim W^\nu_\lambda$ in 
$W^\nu$  at  $\lambda$.
We will show 
that (cf.\ (\ref{eq:rwt2}))
\begin{equation}\label{eq:rres1}
r_\lambda^\nu = \sum_{N\in \cN_\lambda^\nu}
 R(\nu,N), 
\end{equation}
where $R(\nu,N)$ is the integer defined in (\ref{eq:Rdef}).
It will be further shown that the RHS of (\ref{eq:rres1})
now represents
a formal counting  of the 
Bethe vectors of weight $\lambda$
for the \XXZ-type ($U_q(X_n^{(1)})$) spin chain with
 quantum space $W^\nu$ in the $q\rightarrow 0$ limit
(Theorems \ref{thm:C} and \ref{thm:Rnds}).
Therefore, we call the equality (\ref{eq:rres1}) 
the formal  completeness of the \XXZ-type Bethe vectors
(Corollary \ref{cor:completeness3}).
As  the special case of
(\ref{eq:rres1}) for $\nu=(\nu\bk)$, $\nu\bk=\delta_{ab}\delta_{mk}$,
we obtain yet another expression of
$Q\am(x)$ (Theorem \ref{thm:main3}):
\begin{align}\label{eq:qexp2}
Q\am(x)&=\sum_{\lambda\in P}
\left\{\sum_{N\in \cN_\lambda} 
\left(\det_{(b,k),(c,j)} F_{bk,cj}\right)
\prod_{(b,k)}
\frac{1}{N\bk}
\binom{P\bk+N\bk-1}{N\bk-1}
\right\} e^\lambda,\\
F_{bk,cj} &= \delta_{bc}\delta_{kj}P\bk
 + (\alpha_b|\alpha_c)\min(t_c k,t_b j)N_j^{(c)},
\end{align}
where $P\bk$ and $\cN_\lambda$ are in (\ref{eq:pkb1}) and
(\ref{eq:nset5}), and $\det$ and $\prod$ mean the ones
over the index
set $\{\, (b,k)\in H\mid N\bk> 0\, \}$.
The equality (\ref{eq:qexp2}) holds for any $X_n$ 
without assuming the $\cW$ invariance of $Q\am(x)$.

\par
The paper essentially consists of two parts:
In the first part (Sections 2 and 3),
 we derive the number $R(\nu,N)$
in (\ref{eq:rres1}) from a formal counting
 of the \XXZ-type Bethe vectors
 in the $q\rightarrow 0$ limit.
In Theorem \ref{thm:C}
we show that there is a one-to-one correspondence between
 a class of solutions of the Bethe equation
and the one of
the associated linear congruence equation called
the string center equation (SCE).
We then apply the standard M\"obius inversion technique
to count the off-diagonal solutions of the SCE
(Theorem \ref{thm:Rnds}).
In the second part (Sections 4 and 5),
 we show 
the formal completeness of the XXZ-type Bethe 
vectors (\ref{eq:rres1}).
For that purpose we introduce the generating 
series of the numbers $R(\nu,N)$ and derive
its analytic expression (Theorem
\ref{thm:RKgenerating}).
We then show the uniqueness of the series
$\tQ\am$ and that a specialization 
of the above generating series indeed
satisfies the condition in Definition \ref{defn:unique}
(Theorem \ref{thm:main3}).
The two parts are logically almost independent
and therefore able to be read rather separately.


\section{Bethe equation at $q=0$}\label{sec:beq0}

\subsection{The $U_q(X^{(1)}_n)$ Bethe equation}
 \label{subsec:betheeq}

Let $\cN$ be the set in (\ref{eq:ndef1}).
Given $\nu\in \cN$ and a 
sequence of $n$ non-negative integers
$(M_a)_{a=1}^n$,
we associate 
a system of $\sum_{a=1}^n M_a$ equations for 
$\sum_{a=1}^n M_a$ variables
$u_i\sa$ ($a=1,\dots,n; i=1\dots,M_a$),
\begin{equation}\label{eq:bau}
\begin{split}
& \prod_{m=1}^\infty \left(
\frac{\sin\pi\!\!\left(u^{(a)}_i + \frac{\sqrt{-1} m \hbar}{t_a}
\right)}
{\sin\pi\!\!\left(u^{(a)}_i - \frac{\sqrt{-1} m\hbar}{t_a}\right)}
\right)^{\nu^{(a)}_m}\\
=&  - \prod_{b=1}^n\prod_{j=1}^{M_b}
\frac{\sin\pi\!\!\left(u^{(a)}_i - u^{(b)}_j + \sqrt{-1}(\alpha_a
 \vert \alpha_b)\hbar
\right)}
{\sin\pi\!\!\left(u^{(a)}_i - u^{(b)}_j - \sqrt{-1}(\alpha_a \vert
 \alpha_b)\hbar \right)}
\end{split}
\end{equation} 
with $t_a$ and $\alpha_a$ defined in Section \ref{subsec:rect}.
We call (\ref{eq:bau}) the
{\em Bethe (ansatz) equation\/}.
The equation was introduced in \cite{OW,RW}.
It is widely believed  that
for each solution of (\ref{eq:bau}),
one can associate a vector in the space $W^\nu$ 
in (\ref{eq:qspace})
with $X_n$-weight 
\begin{equation}\label{eq:bvw}
 \sum_{(a,m)\in H} m\nu\am \Lambda_a
-\sum_{a=1}^nM_a \alpha_a,
\end{equation}
and the vector is an 
eigenvector (the {\em Bethe vector\/})
of the transfer matrix
of the $U_q(X_n^{(1)})$ spin chain with the quantum space $W^\nu$,
 if it is not the zero vector.

\begin{rem}
Actually, the equation (\ref{eq:bau}) is a special
case of a more general family of the Bethe equations
which depend on the spectral parameters at each site of the
spin chain. Since the analysis below can be easily extended
 for a general case  in a straightforward way,
 we concentrate on the homogeneous case (\ref{eq:bau}).
\end{rem}

By setting $t=\max_{1\leq a \leq n} t_a$,
$q = e^{-2\pi\hbar/t}$,
$x_i\sa=e^{2\pi\sqrt{-1}u_i\sa}$, (\ref{eq:bau}) is  written as 
\begin{equation}\label{eq:ba0}
F_{i+}\sa G_{i-}\sa =-
F_{i-}\sa G_{i+}\sa,
\end{equation}
\begin{alignat}{2}
F_{i+}\sa &=
\prod_{m= 1}^\infty(x_i\sa q^{{mt/t_a}}
-1)^{\nu\am},
&\quad
G_{i+}\sa &=
\prod_{b=1}^n
\prod_{j=1}^{M_b}
(x_i\sa q^{(\alpha_a|\alpha_b)t}
 -x_j\sab),
\\
F_{i-}\sa &=
\prod_{m=1}^\infty(x_i\sa
-q^{mt/t_a})^{\nu\am},
&\quad
G_{i-}\sa&=
\prod_{b=1}^n
\prod_{j=1}^{M_b}
(x_i\sa
 -x_j\sab q^{(\alpha_a|\alpha_b)t}).
\end{alignat}

\subsection{String solution}

We consider  a class of solutions
$(x_i\sa)$ of (\ref{eq:ba0})
such that  $x_i\sa=x_i\sa (q)$ is meromorphic (with respect to
$q$) around
$q=0$.
For a meromorphic function $f(q)$ around $q=0$, 
let $\ord(f)$
 be the order of the leading power of the
Laurent expansion of $f(q)$ around $q=0$, i.e.,
\begin{align}
f(q) &= q^{{\mathrm{ord}}(f)}(f^0 + f^1 q + \cdots ),
 \qquad \ f^0\neq 0,
\end{align}
and let
$\tilde{f}(q) :=  f^0 + f^1q + \cdots$
be the normalized  series.
When $f(q)$ is identically zero, we set $\ord(f)=\infty$.

%
%

\begin{defn}\label{def:admissible}
A meromorphic solution $(x_i\sa)$
of (\ref{eq:ba0}) around $q=0$ is called
 {\em admissible\/} ({\em inadmissible\/}) if
$\ord(F_{i+}\sa G_{i-}\sa) < \infty$ for any $(a,i)$ 
(otherwise).
\end{defn}
For each $N=(N\am)\in \cN$, we set
\begin{align}\label{eq:hdef2}
H'=H'(N):=\{\ (a,m)\in H\mid
N\am > 0\ \},
\end{align}
where $H$ is defined in (\ref{eq:hdef1}).
We have $|H'|<\infty$.
\begin{defn}\label{def:string}\upshape
Let  $(M_a)_{a=1}^n$ be the one in
the Bethe equation (\ref{eq:ba0}),
and let  $N = (N\am)\in \cN$ satisfy
$ \sum_{m=1}^\infty mN\am=M_a$.
A meromorphic solution $(x_i\sa)$ of (\ref{eq:ba0})
around $q=0$ is called a
{\em string solution of pattern
$N$} if \par
(i) $(x_i\sa)$ is admissible.\par
(ii) $(x_i\sa)$ can be arranged as $( x_{m\alpha i}\sa)$ with
\begin{equation}\label{eq:amai1}
 (a,m)\in H', \quad \alpha=1,
\dots, N\am, \quad i=1,\dots,m
\end{equation}
such that\par
(a)  $d_{m\alpha i}\sa
 := {\mathrm{ord}}(x_{m\alpha i}\sa)=(m+1-2i)t/t_a$.
\par
(b) $z^{(a)}_{m \alpha}:=x_{m\alpha1}^{(a)0}=x_{m\alpha2}^{(a)0}=
\cdots = x_{m\alpha m}^{(a)0}\ (\neq 0)$,
where $ x_{m\alpha i}^{(a)0}$ is the coefficient of
the leading power of $x_{m\alpha i}^{(a)}$.
\par\noindent
For each $(a,m,\alpha)$,
$(x_{m\alpha i}\sa)_{i=1}^m$  is called an {\em $m$-string
of color $a$}, and
 $z^{(a)}_{m \alpha}$  is called
the {\em string center} of the $m$-string
$(x_{m\alpha i}\sa)_{i=1}^m$.
(Thus, $N\am$ is the number of the $m$-strings of color $a$.)
\end{defn}
%
%
%
For a string solution
 $(x\sa_{m\alpha i})$,
$x\sa_{m\alpha i}(q)=q^{d\sa_{m\alpha i}}
\tx\sa_{m\alpha i}(q)$, of pattern $N$,
the Bethe equation (\ref{eq:ba0}) reads
\begin{equation}\label{eq:ba2}
F_{m\alpha i+}\sa G_{m\alpha i-}\sa =-
F_{m \alpha i-}\sa G_{m \alpha i+}\sa ,
\end{equation}
\begin{align}
F_{m \alpha i+}\sa&=
\prod_{k=1}^\infty
(\tx_{m \alpha i}\sa q^{d_{m \alpha i}\sa +{k}t/t_a}
-1)^{\nu_k\sa},\\
F_{m \alpha i-}\sa&=
\prod_{k=1}^\infty(\tx_{m \alpha i}\sa q^{d_{m \alpha i}\sa}
-q^{kt/t_a})^{\nu_k\sa},\\
G_{m \alpha i+}\sa&=
\prod_{bk\beta j}
(\tx_{m \alpha i}\sa q^{d_{m \alpha i}\sa +(\alpha_a|\alpha_b)t}
 -\tx_{k \beta j}\sab q^{d_{k\beta j}\sab}),
\\
G_{m\alpha i-}\sa&=
\prod_{bk\beta j}
(\tx_{m \alpha i}\sa q^{d_{m \alpha i}\sa}
 -\tx_{k \beta j}\sab q^{d_{k \beta j}\sab +(\alpha_a|\alpha_b)t}),
\end{align}
where 
$\prod_{bk\beta j}=
\prod_{(b,k)\in H'}
\prod_{\beta = 1}^{N\bk} \prod_{j = 1}^k$,
and the indices
$a$, $m$, $\alpha$, and $i$ run 
in the range (\ref{eq:amai1}).
For a string solution $(x\sa_{m\alpha i})$ of type $N$,
we call (\ref{eq:ba2}) also the Bethe equation.

\par
Notice that $\zeta_{m\alpha i}\sa:=
{\mathrm{ord}}(\tx_{m\alpha i}\sa
 - \tx_{m\alpha i-1}\sa)$
 is positive and finite because of (i) and (ii)
in Definition \ref{def:string}.
We define $y_{m\alpha i}\sa(q)$ ($2 \leq i \leq m$) as
\begin{equation}\label{eq:yzdef1}
q^{{\zeta_{m\alpha i}\sa}}y_{m\alpha i}\sa(q) = \tx_{m\alpha i}\sa(q)
 - \tx_{m\alpha i-1}\sa(q).
\end{equation}
Let us  extract the factors $y_{m\alpha i}\sa$ from
$G_{m\alpha i\pm}\sa$
and introduce $G^{\prime (a)}_{m\alpha i\pm}$ as follows:
\begin{align}
G_{m\alpha i+}\sa &=
\begin{cases}
G^{\prime (a)}_{m\alpha 1+} & i=1\\
G^{\prime (a)}_{m\alpha i+}q^{d_{m\alpha i}\sa
+(\alpha_a|\alpha_a)t+\zeta_{m\alpha i}\sa}y_{m\alpha i}\sa
&2 \le i \le m,
\end{cases}\\
G_{m\alpha i-}\sa &=
\begin{cases}
G^{\prime (a)}_{m\alpha i-}q^{d_{m\alpha i}\sa
+\zeta_{m\alpha i+1}\sa}y_{m\alpha i+1}
&1 \le i \le m-1\\
G^{\prime (a)}_{m\alpha m-} & i=m.
\end{cases}
\end{align}
Now the Bethe equation (\ref{eq:ba2}) takes the form:
\begin{alignat}{3}
\tilde{F}_{m\alpha 1+}\sa
\tilde{G}^{\prime (a)}_{m\alpha 1-}y_{m\alpha 2}\sa
& =
-\tilde{F}_{m\alpha 1-}\sa\tilde{G}^{\prime (a)}_{m\alpha 1+}
&\qquad &i =
1,\label{eq:fgy11}\\
\tilde{F}_{m\alpha i+}\sa
\tilde{G}^{\prime (a)}_{m\alpha i-}y_{m\alpha i+1}\sa
& =
-\tilde{F}_{m\alpha i-}\sa\tilde{G}^{\prime (a)}_{m\alpha i+}
y_{m\alpha i}\sa
&\qquad &2 \le i \le m-1,\label{eq:fgy22}\\
\tilde{F}_{m\alpha m+}\sa\tilde{G}^{\prime (a)}_{m\alpha m-}
& =
-\tilde{F}_{m\alpha m-}\sa
\tilde{G}^{\prime (a)}_{m\alpha m+}y_{m\alpha m}\sa
&\qquad &i = m.\label{eq:fgy33}
\end{alignat}

\subsection{$q\rightarrow 0$ limit of Bethe equation}
\label{subsec:q0be}

Suppose that $(x_{m\alpha i}\sa)$ is a string solution to
 the Bethe equation
(\ref{eq:ba2}).
Taking the leading coefficients 
of (\ref{eq:fgy11})--(\ref{eq:fgy33}),
\begin{alignat}{3}
F^{(a)0}_{m\alpha 1+}
 G^{\prime (a) 0}_{m\alpha 1-}y^{(a)0}_{m\alpha 2}
& =
-F^{(a)0}_{m\alpha 1-}G^{\prime (a) 0}_{m\alpha 1+}
&\qquad &i = 1,\label{eq:fgy1}\\
F^{(a)0}_{m\alpha i+}G^{\prime (a) 0}_{m\alpha i-}
y^{(a)0}_{m\alpha i+1}& =
-F^{(a)0}_{m\alpha i-}G^{\prime (a) 0}_{m\alpha i+}
y^{(a)0}_{m\alpha i}
&\qquad &2 \le i \le m-1,\label{eq:fgy2}\\
F^{(a)0}_{m\alpha m+}G^{\prime (a) 0}_{m\alpha m-} &=
-F^{(a)0}_{m\alpha m-}G^{\prime (a) 0}_{m\alpha m+}
y^{(a)0}_{m\alpha m}
&\qquad &i = m.\label{eq:fgy3}
\end{alignat}
In particular, taking the product of
(\ref{eq:fgy1})--(\ref{eq:fgy3}), we have
\begin{equation}\label{eq:ratio1}
1 = (-1)^m \prod_{i=1}^m
\frac{F^{(a)0}_{m\alpha i+}
G^{(a)0}_{m\alpha i-}}{F^{(a)0}_{m\alpha i-}G^{(a)0}_{m\alpha i+}},
\end{equation}
which turns out to be a key equation.

\subsection{Generic string solution}

We introduce a class of string solutions
for which (\ref{eq:ratio1})
can be explicitly written down.
Let $\xi_{m\alpha i\pm}\sa$,
$\eta_{m\alpha i\pm}\sa$ be the ``superficial'' orders
of the factors $F_{m\alpha i\pm}\sa$,
$G_{m\alpha i\pm}\sa$ in the Bethe equation (\ref{eq:ba2}):
Namely,
\begin{align}
\xi_{m\alpha i+}\sa&=
\frac{t}{t_a}
\sum_{k=1}^\infty
{\nu_k\sa} \min(m+1-2i+k,0),\\
\xi_{m\alpha i-}\sa&=
\frac{t}{t_a}
\sum_{k=1}^\infty
{\nu_k\sa} \min(m+1-2i,k),\\
\eta_{m\alpha i+}\sa &=
\frac{t}{t_a}
\sum_{bk\beta j}
\frac{1}{t_b}\min(t_b(m+1-2i)
+t_at_b(\alpha_a|\alpha_b),t_a(k+1-2j)),\\
\eta_{m\alpha i-}\sa &=
\frac{t}{t_a}
\sum_{bk\beta j}
\frac{1}{t_b}\min(t_b(m+1-2i),t_a(k+1-2j)+t_at_b(\alpha_a|\alpha_b)).
\end{align}
%
%
%
%
\begin{defn}\label{def:generic1}\upshape
A string solution $(x_{m\alpha i}\sa)$ to (\ref{eq:ba2})
is called {\em generic\/} if 
\begin{align}\label{eq:gerdef}
\begin{split}
{\mathrm{ord}}(F_{m\alpha i\pm}\sa) &= \xi_{m\alpha i\pm}\sa,\\
{\mathrm{ord}}(G_{m\alpha i+}\sa) &=
\eta_{m\alpha i+}\sa + \zeta_{m\alpha i}\sa,
\qquad
{\mathrm{ord}}(G_{m\alpha i-}\sa)=
\eta_{m\alpha i-}\sa + \zeta_{m\alpha i+1}\sa,
\end{split}
\end{align}
where $\zeta_{m\alpha1}\sa=\zeta_{m\alpha m+1}\sa=0$.
\end{defn}

Given a quantum space data $\nu
\in \cN$ and a string pattern
$N \in \cN$, we put
\begin{align}
\gamma\am &= \gamma\am(\nu)
 = \sum_{k=1}^\infty \min(m,k) \nu_k\sa,\label{eq:gammadef}\\
P\am &= P\am(\nu,N) = \gamma\am
 - \sum_{(b,k)\in H}
(\alpha_a|\alpha_b) \min(t_b m,t_a k)N\bk.
\label{eq:pdef}
\end{align}
\par
\begin{lem}\label{lem:order1}
We have
\begin{align}\label{eq:ze1}
\begin{split}
& (\xi_{m\alpha i+}\sa+\eta_{m\alpha i-}\sa)
-(\xi_{m\alpha i-}\sa+\eta_{m\alpha i+}\sa)\\
&\qquad =
\begin{cases}
-\frac{t}{t_a}\left(
P_{m+1-2i}\sa+N_{m+1-2i}\sa\right)-\Delta_{m+1-2i}\sa &
1\leq i < \frac{1}{2}(m+1)\\
0 & i = \frac{1}{2}(m+1)\\
\frac{t}{t_a}\left(
P_{2i-m-1}\sa+N_{2i-m-1}\sa\right)+\Delta_{2i-m-1}\sa &
\frac{1}{2}(m+1) < i \leq m,
\end{cases}
\end{split}
\end{align}
where $\Delta_{j}\sa=0$
except for the following cases:
If $t_a=1$ and there is $a'$ such that
$t_{a'}\neq 1$, $C_{aa'}\neq 0$, then
\begin{equation}
\Delta_{j}\sa=
\begin{cases}
-N_{2j}^{(a')}&
t_{a'}=2\\
-
\left(N_{3j-1}^{(a')}
+2N_{3j}^{(a')}
+N_{3j+1}^{(a')}\right)&
t_{a'}=3.
\end{cases}
\end{equation}
\end{lem}
\begin{proof}
It is easy to show them by the case check.
\end{proof}
\par
\begin{prop}\label{pr:order}
A necessary condition for the existence of
a ge\-ner\-ic
 string solution of pattern $N$ is as follows
(\/ $(a,m)\in H'$, $1\le \alpha\leq N\am$, $2 \le i \le m$\/):
\begin{equation}\label{eq:order}
\sum_{k=1}^{\min(i-1,m+1-i)}
\left\{
\frac{t}{t_a}\left(P_{m+1-2k}\sa +
N_{m+1-2k}\sa\right)+\Delta_{m+1-2k}\sa 
\right\}>0.
\end{equation}
\end{prop}
\begin{proof}
Suppose the equation (\ref{eq:ba2})
 admits a generic solution.
Then,
$\xi^{(a)}_{m\alpha i+} +\eta^{(a)}_{m\alpha i-}
 +\zeta^{(a)}_{m\alpha i+1} 
=\xi^{(a)}_{m\alpha i-} +\eta^{(a)}_{m\alpha i+}
 +\zeta^{(a)}_{m\alpha i}$ holds.
Therefore, the LHS of (\ref{eq:ze1}) equals to
$\zeta\sa_{m\alpha i}-\zeta\sa_{m\alpha i+1}$.
Then, the LHS of (\ref{eq:order}) equals to
$\zeta_{m\alpha i}\sa$, which is positive.
\end{proof}

\subsection{String center equation (SCE)}\label{subsec:sce}

Let us compute the equation (\ref{eq:ratio1}) 
for a {\em generic\/} string solution.

\begin{prop}\label{prop:A}
Let $(x_{m\alpha i}\sa)$ be a generic string solution of pattern $N$.
Then its string centers $(z^{(a)}_{m\alpha})$ satisfy
the following equations
$((a,m)\in H', \ 1 \le \alpha \le N\am)$:
\begin{gather}
\prod_{(b,k)\in H'} \prod_{\beta = 1}^{N\bk}
(z^{(b)}_{k\beta})^{A_{am\alpha,bk\beta}} =
 (-1)^{P\am+N\am+1},
\label{eq:sce1}\\
A_{am\alpha,bk\beta} := 
\delta_{ab}\delta_{m k}\delta_{\alpha \beta}(P\am+N\am) +
(\alpha_a|\alpha_b)\min(t_bm,t_ak) - 
\delta_{ab}\delta_{m k}.\label{eq:amat}
\end{gather}
\end{prop}
We call (\ref{eq:sce1}) the {\em string center equation\/}
 (SCE)
of pattern $N$.
The SCE (\ref{eq:sce1})
becomes the linear congruence equation (also called the SCE)
in terms of the variables  $u_{k\beta}$ 
({\it modulo\/} ${\mathbb{Z}}$) defined by
$z^{(b)}_{k\beta} = \exp(2\pi\sqrt{-1}u_{k\beta}\sab)$:
\begin{equation}\label{eq:sce2}
\sum_{(b,k)\in H'}\sum_{\beta=1}^{N\bk}
 A_{am\alpha,bk\beta} u_{k\beta}\sab \equiv
\frac{P\am + N\am + 1}{2} \quad \mathrm{mod}\ {\mathbb{Z}}.
\end{equation}
\begin{proof}
Let us compute the ratio (\ref{eq:ratio1}) explicitly.
\begin{gather*}
\prod_{i=1}^m F_{m\alpha i\epsilon }^{(a)0} =
\begin{cases}
(-1)^{\gamma\am}
\prod_{k=1}^\infty (f^{k}_{am\alpha})^{\nu_k\sa}& \epsilon = +\\
(z^{(a)}_{m\alpha})^{\gamma\am}
\prod_{k= 1}^\infty(f^{k}_{am\alpha})^{\nu_k\sa}& \epsilon = -,
\end{cases}\\
f^{k}_{am\alpha} = \begin{cases}
1 & m \le k\\
(-z^{(a)}_{m\alpha})^{(m-k)/{2}} &
m > k,\  k \equiv m\ \mathrm{mod}\ 2\\
(-z^{(a)}_{m\alpha})^{({m-k-1})/{2}}(z^{(a)}_{m\alpha}-1) &
m > k,\ k \not \equiv m\ \mathrm{mod}\ 2.
\end{cases}
\end{gather*}
In  order to  calculate 
$\prod_{i=1}^m(G^{(a)0}_{m\alpha i-}/G^{(a)0}_{m\alpha
i+})$, it is convenient to evaluate
\begin{gather*}
\begin{split}
&\prod_{i=1}^m\prod_{j=1}^k(\tx_{m\alpha i}\sa
q^{d_{m\alpha i}\sa
+(1+\epsilon)(\alpha_a|\alpha_b)t/2}
 -\tx_{k\beta j}\sab q^{d_{k\beta j}\sab+
(1-\epsilon)(\alpha_a|\alpha_b)t/2})^0\\
&\qquad =
\begin{cases}
(-z^{(b)}_{k\beta})^{
(\alpha_a|\alpha_b)\min(t_b m,t_a k)
-\delta_{a b}\delta_{m k}}g^{bk\beta}_{am\alpha} &
\epsilon = 1\\
(z^{(a)}_{m\alpha})^{(\alpha_a|\alpha_b)
\min(t_b m,t_a k)-\delta_{ab}\delta_{m k}}
(-1)^{(m-1)\delta_{ab}
\delta_{m k}\delta_{\alpha \beta}}g^{bk\beta}_{am\alpha} 
& \epsilon = -1.
\end{cases}\\
\end{split}
\end{gather*}
Here $g_{am\alpha}^{bk\beta}$ are given as follows:
\par
(i) For $t_b(m-1)-t_a(k-1)-t_at_b(\alpha_a|\alpha_b)\equiv 1\ 
\mathrm{mod}\ 2$
\begin{equation*}
g^{bk\beta}_{am\alpha} =
(-z^{(a)}_{m\alpha}z^{(b)}_{k\beta})^{\frac{1}{2}mk-
\frac{1}{2}(\alpha_a|\alpha_b)
\min(t_b m,t_ak)}.
\end{equation*}
\par
(ii) For $t_b(m-1)-t_a(k-1)-t_at_b(\alpha_a|\alpha_b)\equiv 0\
\mathrm{mod}\ 2$,
$(a,m,\alpha)\neq (b,k,\beta)$
\begin{align*}
g^{bk\beta}_{am\alpha} 
&=(-z^{(a)}_{m\alpha}z^{(b)}_{k\beta})^{\frac{1}{2}mk
-\left( \frac{1}{2t_{ab}}
+ \frac{1}{2}(\alpha_a|\alpha_b)\right)
\left(\min(t_b m,t_ak)+ \Delta_{am}^{bk}\right)
+\delta_{a b}\delta_{m k}}\\
&\qquad \times
 (z^{(a)}_{m\alpha}-z^{(b)}_{k\beta})^{\frac{1}{t_{ab}}
\left(\min(t_b m,t_ak)+ (1-t)\Delta_{am}^{bk}\right)
-\delta_{ab}\delta_{m k}},
\end{align*}
where $t_{ab}=\max(t_a,t_b)$, and
$\Delta_{am}^{bk}$ is $0$ except for the cases:
(a) If $t_a < t_b$, $t_bm > t_a k$,
$t_bm-t_a k \equiv \pm 1$ mod $2t$, then
$\Delta_{am}^{bk}$ is $1$ or $-1$ with
$\Delta_{am}^{bk} \equiv t_b m -t_a k$ mod $2t$.
(b) If $t_a > t_b$, $t_bm < t_a k$,
$t_bm-t_a k \equiv \pm 1$ mod $2t$, then
$\Delta_{am}^{bk}$ is $1$ or $-1$ with
$-\Delta_{am}^{bk} \equiv t_b m -t_a k$ mod $2t$.
\par
(iii) For $(a,m,\alpha)=(b,k,\beta)$,
\begin{equation*}
g^{bk\beta}_{am\alpha} =
(-(z^{(a)}_{m\alpha})^2)^{\frac{1}{2}m^2-
\frac{3}{2}m+ 1}
 y^{(a)0}_{m\alpha 2} \cdots y^{(a)0}_{m\alpha m}.
\end{equation*}
The factors $(f^{k}_{am\alpha})^{\nu_k\sa}$ and
$g^{bk\beta}_{am\alpha}$ are
all nonzero for a generic string solution $(x_{m\alpha i}\sa)$.
Thus they are canceled in the ratio in the RHS
of (\ref{eq:ratio1}), and we find
\begin{equation}
\label{eq:ratio}
(-1)^m \prod_{i=1}^m
\frac{F^{(a)0}_{m\alpha i+}G^{(a)0}_{m\alpha i-}}
{F^{(a)0}_{m\alpha i-}G^{(a)0}_{m\alpha i+}}
= (-1)^{P\am +N\am +1}
 \prod_{bk\beta}(z^{(b)}_{k\beta})^{-A_{am\alpha,bk\beta}}.
\end{equation}
{}From (\ref{eq:ratio1}) and (\ref{eq:ratio}) we obtain
(\ref{eq:sce1}).
\end{proof}
%

{}From the conditions
$(f^{k}_{am\alpha})^{\nu_k\sa}\neq 0$
and  $g^{bk\beta}_{am\alpha}\neq 0$
in the above proof,
we see that a string solution is
  generic if and only if its string centers 
$(z_{m\alpha}^{(a)})$ satisfy
the following condition for any $(a,m,\alpha)$:
\begin{equation}\label{eq:nagoya}
\begin{gathered}
\prod_{\scriptstyle k=1\atop
\scriptstyle k \not\equiv m\, (2)}^{m-1}
(z^{(a)}_{m\alpha} - 1 )^{\nu_k\sa} \neq 0,\\
\prod_{\scriptstyle 
bk \beta  \, (\neq am\alpha)\atop
{\scriptstyle
t_b(m-1)-t_a(k-1)
\atop
\scriptstyle 
\qquad -t_at_b (\alpha_a|\alpha_b)\equiv 0\, (2)}}
(z^{(a)}_{m\alpha} - z^{(b)}_{k\beta})^{\frac{1}{t_{ab}}
\left(\min(t_b m,t_ak)+ (1-t)\Delta_{am}^{bk}\right)
-\delta_{ab}\delta_{m k}} \neq 0.
\end{gathered}
\end{equation}

\begin{defn}
A solution to the SCE (\ref{eq:sce1}) is called {\em generic\/}
if it satisfies the  condition (\ref{eq:nagoya}).
\end{defn}

Let $A$ be the matrix 
with  the entry $A_{am\alpha, bk\beta}$ in (\ref{eq:amat}).

\begin{prop}\label{prop:B}
Suppose that $N \in \cN$ satisfies the conditions (\ref{eq:order})\/
and\/ $\det A\neq 0$.
Then, for each generic solution $(z^{(a)}_{m\alpha})$ to
the SCE of pattern $N$,
there exists a unique generic string solution
 $({x}^{\prime (a)}_{m\alpha i}(q))$ of pattern $N$ to the
Bethe equation (\ref{eq:ba2}) such that
its string center $z^{\prime (a)}_{m\alpha}$
is equal to $z^{(a)}_{m\alpha}$.
\end{prop}

To prove Proposition \ref{prop:B},
we introduce new  variables
 $w_{m\alpha i}\sa $ as
\begin{equation}
w_{m\alpha i}\sa =
\begin{cases}
\tx_{m\alpha i}\sa & i=1\\
y_{m\alpha i}\sa & 2 \le i \le m.
\end{cases}
\end{equation}
Then 
\begin{equation}\label{eq:wxyz}
\tx_{m\alpha i}\sa = w_{m\alpha 1}\sa
 + q^{\zeta_{m\alpha 2}\sa}w_{m\alpha 2}\sa + \cdots
+ q^{\zeta_{m\alpha i}\sa}w_{m\alpha i}\sa \qquad 1 \le i \le m.
\end{equation}
Let us write the $i$th equation of (\ref{eq:fgy11})--(\ref{eq:fgy33})
as $L_{m\alpha i}\sa = R_{m\alpha i}\sa$.
Let $J = (J_{am\alpha i, bk\beta j})$ be a matrix 
with  entry
$J_{am\alpha i, bk\beta j}
 = \frac{\partial}{\partial w_{k\beta j}\sab}
\bigl(\frac{L_{m\alpha i}\sa}{R_{m\alpha i}\sa}-1\bigr)$.
\begin{lem}\label{lem:jacobian}
If $N\in \cN$ satisfies the conditions (\ref{eq:order})\/
and\/ $\det A\neq 0$,
then  $\det J $ is not zero at $q=0$.
\end{lem}
\begin{proof}
Owing to the assumption (\ref{eq:order}), we have
$\zeta_{m\alpha i}\sa  > 0$ for any $(m,\alpha,i)$.
Since $L^{(a)0}_{m\alpha i} = R^{(a)0}_{m\alpha i} \neq 0$,
it suffices to show $\det{\cJ} \neq 0$ for 
${\cJ}_{am\alpha i, bk\beta j} =
\frac{\partial}{\partial w_{k\beta j}\sab}\log
\frac{L_{m\alpha i}\sa}{R_{m\alpha i}\sa}
$.
{}From (\ref{eq:wxyz}) both
$\frac{\partial\tilde{F}_{m\alpha i \pm}\sa}
{\partial w_{k\beta j}\sab}$ and
$\frac{\partial\tilde{G}^{\prime(a)}_{m\alpha i \pm}}
{\partial w_{k\beta j}\sab}$
for  $j \neq 1$ are zero at $q=0$.
Thus among ${\cJ}^0_{am\alpha i, bk\beta j}$'s
 the non-vanishing
 ones are only
${\cJ}^0_{am\alpha i, bk\beta 1}$ $(1\leq i \leq m)$,
${\cJ}^0_{am\alpha i, am\alpha i}
 = -1/{y^{(a)0}_{m\alpha i}}$ $(2 \le i \le m)$,
and 
${\cJ}^0_{am\alpha i, am\alpha i+1}
= 1/{y^{(a)0}_{m\alpha i+1}}$ $(1 \le i \le m-1)$.
Let $\vec{{\cJ}}^{0}_{am\alpha i}
 = ({\cJ}^0_{am\alpha i,
bk\beta j})_{bk\beta j}$
be the $(am\alpha i)$-th row vector of the matrix ${\cJ}$.
In view of the above result,
the linear dependence $\sum_{am\alpha i}c_{am\alpha
i}\vec{{\cJ}}^{0}_{am\alpha i} = 0$
can possibly hold only when $c_{am\alpha i}$ is independent of $i$.
Consequently we consider the equation
$\sum_{am\alpha}
c_{am\alpha}\sum_{i=1}^m \vec{{\cJ}}^{0}_{am\alpha
i} = 0$.
The $(bk\beta 1)$-th component of the vector
$\sum_{i=1}^m \vec{{\cJ}}^{0}_{am\alpha i}$
is given by
\begin{displaymath}
\lim_{q \rightarrow 0} \frac{\partial}{\partial w_{k\beta 1}\sab}
 \log
\prod_{i=1}^m \frac{-\tilde{F}_{m\alpha i+}\sa
\tilde{G}^{\prime(a)}_{m\alpha i-}}
{\tilde{F}_{m\alpha i-}\sa\tilde{G}^{\prime(a)}_{m\alpha i+}}
= \frac{\partial}{\partial z^{(b)}_{k\beta}} \log
\prod_{i=1}^m \frac{-F^{(a)0}_{m\alpha i+}G^{(a)0}_{m\alpha i-}}
{F^{(a)0}_{m\alpha i-}G^{(a)0}_{m\alpha i+}},
\end{displaymath}
where we have taken into account (\ref{eq:wxyz}) and
$\zeta_{m\alpha i}\sa > 0$.
Due to (\ref{eq:ratio}) the last expression is equal to
$-A_{am\alpha,b k\beta}/z^{(b)}_{k\beta}$.
Therefore the equation
$\sum_{am\alpha}c_{am\alpha}
 \sum_{i=1}^m \vec{{\cJ}}^{0}_{am\alpha i} = 0$
is equivalent to
$\sum_{am\alpha}c_{am\alpha}
A_{am\alpha,b k\beta} = 0$ for any $(b,k,\beta)$.
This admits only the trivial solution for
$c_{am\alpha}$ if $\det A \neq 0$.
\end{proof}

{\em Proof of Proposition \ref{prop:B}}.
The Bethe equations (\ref{eq:fgy11})--(\ref{eq:fgy33})
are simultaneous equations in the variables
$((w_{m\alpha i}\sa), q)$.
At $q=0$, (\ref{eq:fgy11})--(\ref{eq:fgy33})
reduce to (\ref{eq:fgy1})--(\ref{eq:fgy3}).
The latter fix $(y^{(a)0}_{m\alpha i})$ unambiguously
once a generic solution $(z^{(a)}_{m\alpha})$ to the 
SCE is given.
Denote the resulting value of $w_{m\alpha i}\sa$ by
$w^{(a)0}_{m\alpha i}$.
{}From Lemma \ref{lem:jacobian} and
the implicit function theorem, there uniquely
exist the functions $w^{\prime (a)}_{m\alpha i}(q)$ satisfying
(\ref{eq:fgy11})--(\ref{eq:fgy33}) and
$w^{\prime (a)0}_{m\alpha i} = w^{(a)0}_{m\alpha i}$.
\qed
\par
{} From Propositions \ref{prop:A} and \ref{prop:B},
we obtain the main statement in this section.
\begin{thm}\label{thm:C}
Suppose that $N \in \cN$ satisfies the conditions (\ref{eq:order})\/
and\/ $\det A\neq 0$.
Then, there is a one-to-one correspondence between
generic string solutions of pattern $N$ 
to the Bethe equation (\ref{eq:ba2})
and generic solutions to the SCE (\ref{eq:sce1}) of pattern $N$.
\end{thm}

\begin{rem}\label{rem:order}
As we will see later in Lemma \ref{lem:positive},
the condition $\det A\neq 0$ in Theorem \ref{thm:C} is satisfied if
$N\in \cN$ satisfy the condition:
\begin{equation}\label{eq:order2}
\text{
$P\am(\nu,N) \ge 0$ for any
$(a,m)\in H'$.}
\end{equation}
More strongly, $\det A$ is positive under (\ref{eq:order2}).
In general, the conditions (\ref{eq:order})
and (\ref{eq:order2})  are simultaneously
satisfied if
$\sum_m mN\am$ is sufficiently smaller than
$\sum_m \nu\am$.
Because of (\ref{eq:bvw}),
$\sum_{(a,m)\in H} mN\am \alpha_a$ measures the difference between the
weight of the corresponding Bethe vector and the
highest weight  of the
quantum space $W^\nu$.
Thus,
if $\sum_m \nu\am$ are large enough, 
the conditions in Theorem \ref{thm:C} are satisfied
at least ``near the highest weight''.
\end{rem}


\section{Counting of off-diagonal solutions to SCE}
\label{sec:counting}

In this section the off-diagonal
solutions of the SCE will be counted
under a certain condition (Theorem \ref{thm:Rnds}).

\subsection{Off-diagonal solution}\label{subsec:offdiag}

In what follows, the symbol $\binom{k}{j}$
($k\in \bC$, $j\in \bZ$) will 
denote the binomial coefficient:
\begin{equation}
\binom{k}{j} = \begin{cases}
k(k-1) \cdots (k-j+1)/j! &  j >0 \\
1 &  j = 0\\
0 & j <0.
\end{cases}
\end{equation}

For each $\nu$, $N\in \cN$, we define the number
$R(\nu,N)$ as follows:
For $N\neq 0\in \cN$, 
\begin{align}
R(\nu,N) &=
\left(\det_{(a,m),(b,k) \in {H'}}F_{am,bk}\right)
\prod_{(a,m) \in H'} \frac{1}{N\am}
\binom{P\am + N\am - 1}{N\am - 1},
\label{eq:Rdef}\\
\label{eq:fdef}
F_{am,bk} &= \delta_{ab}\delta_{mk}P\am
 + (\alpha_a|\alpha_b)\min(t_b m,t_a k)N\bk,
\end{align}
where $H'=H'(N)$ and $P\am = P\am(\nu,N)$
are given by (\ref{eq:hdef2}) and (\ref{eq:pdef}).
For $N=0$, 
we set $R(\nu,0) = 1$ irrespective of $\nu$.
For any  $N\in \cN$,
$R(\nu,N)$
is  an integer (cf.\ Lemma 
\ref{lem:ralt}), though it is not always a positive one.

\begin{defn}\label{def:off2}\upshape
A solution $(z_{m\alpha}^{(a)})$ to the SCE
is called {\em off-diagonal\/} ({\em diagonal\/}) if
$z^{(a)}_{m\alpha} = z^{(a)}_{m\beta}$ only for $\alpha=\beta$
(otherwise).
\end{defn}
Notice that ``off-diagonal'' above is weaker than
the condition that $z^{(a)}_{m\alpha}$'s are {\em all distinct}.

\begin{thm}\label{thm:Rnds}
If $N\, (\neq 0)\in \cN$
satisfies the condition (\ref{eq:order2}),
then the number of off-diagonal solutions to
the SCE (\ref{eq:sce1}) of pattern $N$
divided by $\prod_{(a,m)\in H'} N\am!$ is equal to $R(\nu,N)$.
\end{thm}

\begin{rem}
In contrast with Definition \ref{def:off2},
let us call a solution to the Bethe equation (\ref{eq:ba2})
{\em off-diagonal\/} if $x^{(a)}_{m\alpha i}(q)$'s
 are all distinct.
Theorem \ref{thm:Rnds} is motivated by the fact
for $X_n=A_1$ \cite{TV} that 
(i) the Bethe vector  associated
to each solution to the Bethe equation 
does not vanish if and only if
 $x_i^{(1)}$'s are admissible and off-diagonal;
(ii) the Bethe vector is invariant under the
permutations of $x_i^{(1)}$'s.
Combining Theorem \ref{thm:Rnds} with Theorem \ref{thm:C},
the number $R(\nu,N)$ correctly counts
the off-diagonal string solutions of
pattern $N$ (modulo permutation) to the Bethe
equation (\ref{eq:ba2})
if $N$ satisfies all
the following conditions:
\par
(i) $N\in \cN$ satisfies the conditions (\ref{eq:order})
and (\ref{eq:order2}).
\par
(ii) All the off-diagonal string solutions 
of pattern $N$ to the Bethe equation  are generic.
\par
(iii) For each off-diagonal string solution
of pattern $N$ to the Bethe equation,
the string centers are all distinct.
\par
(iv) For each off-diagonal solution to the SCE of pattern $N$,
$z_{m\alpha}^{(a)}$'s are all distinct.
\par\noindent
Unfortunately, so far we do not  know a more general condition 
where the one-to-one correspondence between
the off-diagonal string solutions to the Bethe equation
(\ref{eq:ba2})
and the off-diagonal solutions to the SCE (\ref{eq:sce1})
holds.
\end{rem}

\subsection{Proof of Theorem \ref{thm:Rnds}}
\label{subsec:pr32}

We work with the logarithmic form of the SCE (\ref{eq:sce2}),
\begin{equation}\label{eq:sce3}
A \vec{u} \equiv \vec{c} \quad \mathrm{mod}\ \mathbb{Z}^d,
\end{equation}
where $d = \sum_{(a,m)\in H'}N\am$.
Here $\vec{u} = (u_{k\beta}\sab)$ is the unknown and $\vec{c}$
is some constant vector.
The matrix $A = (A_{am\alpha, bk\beta})$ is
specified by (\ref{eq:amat}).
We consider the solutions  $u_{k \beta}\sab$ {\em modulo\/}
integers.

The following fact is well-known (cf.\ \cite[1.2.2, Lemma 1]{C}):
\begin{lem}\label{lem:basic}
Let $B$ be an $r$ by $r$ integer matrix with $\det B \neq 0$.
Then for any $\vec{b} \in \mathbb{R}^r$, the equation
$B \vec{x} \equiv \vec{b}$ $\mathrm{mod}$ $\mathbb{Z}^r$ has exactly 
$\vert \det B \vert$ solutions $\vec{x}$ 
in $(\mathbb{R}/\mathbb{Z})^r$.
\end{lem}

Therefore, if $\det A\neq 0$, then
the number of the solutions to (\ref{eq:sce3}), including
diagonal ones, is given by $|\det A|$.
One can systematically remove the diagonal solutions
using the M\"obius inversion method \cite{A,S} as follows.

For a given positive integer $K$,
$\pi = (\pi_1, \ldots, \pi_l)$ is called a {\em partition\/} of
a set $\{\,1, \ldots, K\,\}$ if
\begin{equation*}
\{\, 1, \ldots, K\,\} = \pi_1 \sqcup \cdots \sqcup \pi_l
\end{equation*}
is a disjoint union decomposition.
Here, the ordering of $\pi_1$, $\pi_2$, \dots, $\pi_l$ is
ignored.
Each $\pi_i$ is called a {\em block\/} of $\pi$,
and $l$ is called a {\em length\/} of $\pi$.
Let $L_K$ denote the set of partitions of $\{\,1, \ldots, K\,\}$.
$L_K$ becomes a partially ordered set (poset)
by the following partial order:
Given two partitions $\pi, \pi' \in L_K$, we say
$\pi \le \pi'$ if each block of $\pi'$ is 
contained in a block of $\pi$.
Let $\mu(\pi,\pi')$ be the M\"obius function for the poset $L_K$.
It is well-known \cite{A,S}  that
\begin{lem}\label{pr:inversion1}
Let $X$ be an indeterminate. For any $\pi \in L_K$ we have
\begin{align}
( X )_{l(\pi)} &= \sum_{\pi' \le \pi} \mu(\pi',\pi) X^{l(\pi')},
\end{align}
where $( X )_l = X(X-1) \cdots (X-l+1)$.
\end{lem}

For a given string pattern  $N\in \cN$,
we consider the direct product of posets
$\cL_N=\prod_{(a,m)\in H'} L_{N\am}$;
for $\pi=(\pi\am)$, 
$\pi'=(\pi^{\prime(a)}_m)\in \cL_N$,
we define $\pi \leq \pi'$ when
$\pi\am \leq \pi^{\prime(a)}_m$ for each $(a,m)$.
Below, $\mu(\pi,\pi')$ means
the M\"obius function for $\cL_N$.
We set $l(\pi)=\sum_{(a,m)\in H'} l(\pi\am)$.
For each $\pi=(\pi\am)\in \cL_N$, let
\begin{displaymath}
\begin{split}
\Sol'_\pi = \{\,\vec{u} = (u_{k\beta}\sab)
 \mid &\, \text{$\vec{u}$ is a solution of (\ref{eq:sce3})},\
u_{k\alpha}\sab = u_{k\beta}\sab
\text{ if}\\
& \text{ $\alpha$ and $\beta$ belong 
to the same block of $\pi\bk$}\,\},\\
\Sol_\pi = \{\,\vec{u} = (u_{k\beta}\sab)
 \mid & \, \text{$\vec{u}$ is a solution of (\ref{eq:sce3})},\
 u_{k\alpha}\sab = u_{k\beta}\sab
 \text{ if and only if}\\
&\text{ $\alpha$ and $\beta$ belong to the same
block of $\pi\bk$}\,\}.
\end{split}
\end{displaymath}
In particular, 
$ \Sol_{\pi_{\text{max}}}$
is the set of the off-diagonal solutions of (\ref{eq:sce3}),
where $\pi_{\text{max}}$
is the maximal element in $\cL_N$.
\begin{lem}
\begin{equation}\label{eq:sol1}
\vert \Sol_\pi \vert = \sum_{\pi' \le \pi} \mu(\pi',\pi)
\vert \Sol'_{\pi'} \vert\qquad \end{equation}
\end{lem}
\begin{proof}
By definition
$\vert \Sol'_\pi \vert = \sum_{\pi' \le \pi} \vert
 \Sol_{\pi'} \vert.$
Applying the M\"obius inversion formula \cite{S}, we obtain
(\ref{eq:sol1}).
\end{proof}
\par
To use the formula (\ref{eq:sol1}),
let us  evaluate $\vert \Sol'_\pi \vert$.
With the constraint,
$u_{k\alpha}\sab = u_{k\beta}\sab$ on $\vec{u} = (u_{k\beta}\sab)$
 if
$\alpha$ and $\beta$  belong to the same block of
$\pi\bk$, 
the SCE (\ref{eq:sce3}) reduces to the following form:
\begin{equation}
A^\pi \vec{u}_\pi \equiv \vec{c}_\pi \quad 
\mathrm{mod}\
\mathbb{Z}^{l(\pi)}
\label{eq:api}
\end{equation}
In the new unknown  $\vec{u}_\pi = (u_{k\beta}\sab)$,
$\beta$ is now labeled by the blocks of $\pi\bk$.
The matrix
$A^\pi$ is an  integer matrix of size $l(\pi)$  obtained by
a reduction of the matrix $A$ as follows:
It is formed by summing up the $(bk\beta)$-th columns of $A$
over those $\beta$ belonging to the same block of $\pi^{(b)}_k$,
and discarding all but one rows for each  block.
Explicitly,
\begin{equation}
A^\pi_{ami, bkj} = \delta_{ab}\delta_{mk}\delta_{ij}(P\am+N\am) +
\{(\alpha_a|\alpha_b)\min(t_b m,t_a k)
-\delta_{ab}\delta_{mk}\}\vert \pi^{(b)}_{k,j} \vert
\end{equation}
where $1 \le i \le l(\pi_m^{(a)})$ and
$1 \le j \le l(\pi_k^{(b)})$.
{}From Lemma \ref{lem:basic} we have
$\vert \Sol'_\pi \vert = \vert \det A^\pi \vert$ if
$\det A^\pi \neq 0$.
\par
It is easy to show the following formula by 
elementary transformations of $A^\pi$.
\begin{lem}\label{lem:detapi}
\begin{align}
\det A^\pi &= \left(\det_{(a,m),(b,k) \in H'}
(F_{am,bk})\right)
\prod_{(a,m) \in H'} (P\am
 + N\am)^{l(\pi\am)-1}.\label{eq:detapi}
\end{align}
\end{lem}
Furthermore
\begin{lem}\label{lem:positive}
If $N\, (\neq 0)\in \cN$
satisfies the condition (\ref{eq:order2}),
then\/ $\det A^\pi > 0$.
\end{lem}
\begin{proof}
By {Lemma} \ref{lem:detapi}, it suffices to verify
$\det_{H'}F > 0$.
Actually, a stronger statement holds: 
Let us forget the relation (\ref{eq:pdef}),
and regard $P\am$  in (\ref{eq:fdef})
as a {\em nonnegative\/} integer which is independent of $N\am$'s.
Then, $\det_{H'}F > 0$ still holds.
We prove the last statement by the double induction
on $\vert H' \vert$ and the sum
$\sum_{(a,m)\in H'} P\am$.
First, let $H'$  ($\neq \emptyset$) be arbitrary,
and suppose
$\sum_{(a,m)\in H'} P\am=0$ (i.e., $P\am=0$ for any $(a,m)
\in H'$).
Then $\det_{H'}F > 0$ is equivalent to
the positivity of 
$\det_{H'} (\alpha_a|\alpha_b)\min(m/t_a, k/t_b)$,
which is a principal minor of the tensor product of two
positive-definite 
matrices, $((\alpha_a|\alpha_b))_{1\leq a,b\leq n}$ and
$(\min(m,k)/s)_{1\leq m,k\leq L}$
with $L$ and $s$ some integers.
Therefore, the claim is true.
Next, suppose that 
$H'=\{ (a,m) \}$ (i.e., $|H'| =1$)
and $P\am $ is any nonnegative integer.
Then $\det_{H'}F = P\am + (\alpha_a|\alpha_a)
m t_a N\am > 0$.
Finally, let $H'$ ($|H'|\geq 2$) and
 $\sum_{(a,m)\in H'}P\am$ ($ > 0)$ be arbitrary.
Then there exists some
$(r,i) \in H'$ such that $P_i^{(r)} > 0$.
Set $\tilde{P}^{(a)}_m = P\am - 1$ if $(a,m)=(r,i)$,
 $\tilde{P}_m^{(a)}=P\am$ otherwise,
and
$\tilde{H}' = H' \setminus  \{(r,i)\}$.
Then, one can split the determinant as
$\det_{H'}F(\{P\am \}) =
\det_{\tilde{H}'}F(\{P\am \}) + \det_{H'}F(\{\tilde{P}^{(a)}
_m \})$.
By the induction hypothesis, the RHS is positive.
\end{proof}
\par

Assembling Lemmas \ref{lem:basic}--\ref{lem:positive}, we have
\begin{align*}
\frac{|\Sol_{\pi_{\text{max}}}|}{\prod_{(a,m)\in H'}
N\am!}&=
\frac{1}{\prod_{(a,m) \in H'}N\am!}
\sum_{\pi \in \cL_N}
\mu(\pi,\pi_{\text{max}})
\det A^\pi \\
&=
\left(\det_{(a,m),(b,k) \in H'}(F_{am,bk})\right)
\prod_{(a,m) \in H'}
\frac{(P\am+N\am)_{N\am}}{N\am!(P\am+N\am)}\\
&=  R(\nu,N).
\end{align*}
This completes the proof of Theorem \ref{thm:Rnds}.


\section{Generating series}\label{sec:multiplicity}

The number $R(\nu,N)$, which was introduced 
in Section \ref{sec:counting} to count
the \XXZ-type Bethe vectors,
is our main 
concern below.
Let us recall the  definition of 
$R(\nu,N)$ in (\ref{eq:Rdef}),
where $\nu$, $N\in \cN$, and $H$
and $\cN$ are defined in (\ref{eq:hdef1}) and
(\ref{eq:ndef1}):
For $N\neq 0$,
\begin{align}
\label{eq:Rdef2}
R(\nu,N) &=
\left(
\det_{(a,m),(b,k) \in H'}F_{am,bk}\right)
\prod_{(a,m) \in H'} \frac{1}{N\am}
\binom{P\am + N\am - 1}{N\am - 1},\\
H' & = H'(N)=\{\ (a,m)\in H\mid N\am >0\ \},\\
\gamma\am &= \gamma\am(\nu)
 = \sum_{k=1}^\infty \min(m,k) \nu_k\sa,
\label{eq:gammadef2}\\
\label{eq:pdef2}
P\am &= P\am(\nu,N) = \gamma\am
 - \sum_{(b,k)\in H}
(\alpha_a|\alpha_b) \min(t_b m,t_a k)N\bk,\\
\label{eq:fdef2}
F_{am,bk} &= \delta_{ab}\delta_{mk}P\am
 + (\alpha_a|\alpha_b)\min(t_b m,t_a k)N\bk.
\end{align}
For $N=0$, we set $R(\nu,0)=1$.
{}From now on, we forget the conditions
(\ref{eq:order}) and (\ref{eq:order2})
 required for 
 Theorems \ref{thm:C} and \ref{thm:Rnds}.
For any  $N\in \cN$,
$R(\nu,N)$
is an integer (cf.\ Lemma 
\ref{lem:ralt}), though it is not always a positive one.
\par 
On the other hand, the number
\begin{equation}\label{eq:Kdef}
K(\nu,N) = \prod_{(a,m)\in H}
\binom{P\am + N\am}{N\am}
\end{equation}
was  introduced in \cite{KR}
to count the \XXX-type Bethe vectors.
Below we treat two numbers, $R(\nu,N)$ and $K(\nu,N)$, in a parallel
way so that the relation between them becomes transparent.

\subsection{Generating series}
It is natural to introduce
the {\em generating series\/} of
$R(\nu,N)$ and $K(\nu,N)$,
\begin{align}
\cR^\nu(w) &= \sum_{N\in \cN} R(\nu,N) w^N,
\qquad
w^N=\prod_{(a,m)\in H} (w\am)^{N\am},
\label{eq:rwdef}\\
\cK^\nu(w) &= \sum_{N\in \cN} K(\nu,N)w^N.\label{eq:kwdef}
\end{align}
It turns out that the following truncation (projection)
is appropriate
for our purpose (cf.\ Remark \ref{rem:main}):
Let $l$ be  a fixed nonnegative integer $l$, and let
\begin{align}
\label{eq:hl1}
H_l&=\{\, (a,m) \mid 1\leq a\leq n,\ 1\leq m \leq t_a l\,\},\\
\label{eq:nl1}
\cN_l&=
\{\, N=(N\am)_{(a,m)\in H}\mid
N\am \in \mathbb{Z}_{\ge 0},\
N\am = 0\ \text{for}\ (a,m)\notin H_l\,\},
\end{align}
so that $\varinjlim H_l=H$ and $\varinjlim \cN_l=\cN$.
The truncated generating series are defined as 
\begin{align}
\cR^\nu_l(w) &= \sum_{N\in \cN_l} R(\nu,N) w^N,
\qquad
w^N=\prod_{(a,m)\in H_l} (w\am)^{N\am},
\label{eq:rlwdef}\\
\cK^\nu_l(w) &= \sum_{N\in \cN_l} K(\nu,N)w^N.
\label{eq:klwdef}
\end{align}

\begin{prop}\label{prop:nuzero}
\begin{equation*}
\cR^0_l(w) = 1.
\end{equation*}
\end{prop}
\begin{proof}
By definition, $R(0,0)=1$.
Let $\nu=0$ and $N\neq 0$.
{}From (\ref{eq:fdef2}), we have
$\sum_{(b,k)\in H'}F_{am,bk} = \gamma\am$,
and $\gamma\am = 0$ when $\nu = 0$.
Thus, $\det F=0$. 
Therefore,  $R(0, N) = 0$ by (\ref{eq:Rdef2}).
\end{proof}
\par
In contrast, $\cK^0_l(w)$ is not so simple.
See (\ref{eq:K2}) and (\ref{eq:R2}).

We need the following alternative expression of $R(\nu,N)$.

\begin{lem}\label{lem:ralt}
For any $\nu$, $N\in {\cN}$, the following equality holds:
\begin{align}\label{eq:Rldef3}
\begin{split}
 R(\nu,N) 
&= \sum_{J \subset H}
\left\{ D_J
\prod_{(a,m) \in H}\binom{P[J]\am + N[J]\am}{N[J]\am}
\right\},
\end{split}
\end{align}
where the sum is taken over all the finite subsets $J$ of $H$,
\begin{align}
D_J &= \begin{cases}
1 & \text{ if } J = \emptyset\\
\det_{(a,m),(b,k) \in J}
D_{am,bk}
 & \text{otherwise},
\end{cases}\label{eq:Ddef}\\
\label{eq:Dmat}
D_{am,bk}
&=
(\alpha_a|\alpha_b)\min(t_b m,t_a k)-
\delta_{ab}\delta_{mk},
\end{align}
$P[J]\am:=P\am(\nu[J],N[J])$ with
\begin{alignat}{2}
\label{eq:nuJ}
\nu[J] &= (\nu[J]\am),\qquad & 
\nu[J]\am &= \nu\am
 - \sum_{(b,k)\in J}(\alpha_a|\alpha_b)
B_{bk,am},\\
N[J] &= (N[J]\am),\qquad &
N[J]\am &= N\am - \theta((a,m) \in J),
\label{eq:nJ}
\end{alignat}
$B_{bk,am}$ is defined in (\ref{eq:bfunc}),
and\/
$\theta(\text{\rm true}) = 1$, 
$\theta(\text{\rm false}) = 0$.
\end{lem}
\begin{proof}
If $N=0$, then the both hand sides of 
(\ref{eq:Rldef3}) is 1.
Suppose that $N\neq 0$.
By (\ref{eq:fdef2}) and (\ref{eq:Dmat}),
we have 
$F_{am,bk}=\delta_{ab}\delta_{bk}(P\am+N\am)+D_{am,bk}
N\bk$.
By
splitting the sum $(P\am+N\am) + D_{am,am}N\am$
of each diagonal element
in $\det F$ in (\ref{eq:Rdef2}), 
$R(\nu,N)$ is written as
\begin{align}
\label{eq:dj1}
 \sum_{J \subset H'}
\left\{ D_J
\prod_{(a,m) \in H'\setminus J}\binom{P\am + N\am}{N\am}
\prod_{(a,m) \in J}\binom{P\am + N\am - 1}{N\am - 1}\right\}.
\end{align}
Then, (\ref{eq:Rldef3}) follows from (\ref{eq:dj1})  
and the fact
\begin{equation}\label{eq:peq}
\begin{split}
P[J]\am- P\am&=
-\sum_{j=1}^{\infty} \sum_{(b,k)\in J}
\min(m,j)(\alpha_a |\alpha_b) B_{bk,aj}\\
&\qquad +
\sum_{(b,k)\in H}
(\alpha_a | \alpha_b)\min(t_bm,t_a k)\theta((b,k)\in J)=0,
\end{split}
\end{equation}
where the last equality in (\ref{eq:peq}) is due to (\ref{eq:bsum}).
\end{proof}

\begin{rem}\label{rem:Next}
Let
\begin{align}
\label{eq:tcN}
\begin{split}
\tilde{\cN}=\{\, N=(N\am)_{(a,m)\in H}
\mid\, &\text{$N\am\in \bZ$, $N\am=0$ except for}\\
& \,\text{finitely-many $(a,m)$}\, \}.
\end{split}
\end{align}
We extend the definition of $R(\nu,N)$ to $N\in \tilde{\cN}$
such that $R(\nu,N)=0$ if $N\in \tilde{\cN} \setminus \cN$.
Then, the equality (\ref{eq:Rldef3}) still
holds for any $N\in \tilde{\cN}$
because the both hand sides of 
(\ref{eq:Rldef3}) is 0 if $N\in \tilde{\cN}\setminus \cN$.
We use this fact in the proof of Lemma \ref{lem:sum}.
\end{rem}

For $N\in \cN_l$,
the expression (\ref{eq:Rldef3}) reads as
\par
\begin{lem}\label{lem:ralt2}
For $N\in \cN_l$, we have
\begin{align}
R(\nu,N) = \sum_{J \subset H_l}
\left\{ D_J
\prod_{(a,m) \in H_l}\binom{P[J]\am + N[J]\am}{N[J]\am}
\right\}.\label{eq:Rldef4}
\end{align}
\end{lem}

\begin{prop}\label{prop:RKconv}
$\cR^\nu_l(w)$ and $\cK^\nu_l(w)$ converge
for $|w\am|<{(2m-1)^{2m-1}}  / {(2m)^{2m}}$.
\end{prop}

\begin{proof}
First, we consider $\cK^\nu_l(w)$.
For given $(a,m)\in H_l$ and $N\in \cN_l$,
let $N'\in \cN_l$ be $N_k^{\prime(b)}
=N\bk + \delta_{ab}\delta_{mk}$,
and $P_k^{\prime (b)}=P\bk(\nu,N')$.
Then, it is easy to check that
\begin{align*}
\lim_{N\am\to \infty}
\binom{P\bk+N\bk}{N\bk}\Bigm/
\binom{P_k^{\prime (b)}+N_k^{\prime(b)}}
{N_k^{\prime(b)}}
=
\begin{cases}
-\frac{(2m-1)^{2m-1}}{(2m)^{2m}}&  (b,k)=(a,m)\\
1 &  (b,k)\neq (a,m).
\end{cases}
\end{align*}
Therefore, $\cK^\nu_l(w)$ converges for
$|w\am|<{(2m-1)^{2m-1}}/{(2m)^{2m}}$.
Next, consider $\cR^\nu_l(w)$.
By Lemma \ref{lem:ralt2},
$\cR^\nu_l(w)$ is a linear sum of the
power series  whose coefficient
of $w^N$ is 
$\prod_{(a,m)\in H_l}\binom{P[J]\am+N[J]\am}{N[J]\am}$.
Again, each series converges  for
$|w\am|<{(2m-1)^{2m-1}}/{(2m)^{2m}}$.
\end{proof}

\subsection{Basic identity}\label{subsec:basic}

\par
To proceed, we use an identity found
by \cite{K2} for $A_n$
and generalized to $X_n$ by \cite{HKOTY}.

\par

Let 
 $v= (v\am)_{(a,m)\in H_l}$ and
$z = (z\am)_{(a,m)\in H_l}$ 
be complex multivariables.
We define a map $z=z(v)$ as
\begin{align}\label{eq:zvrel}
z\am (v) &= v\am
 \prod_{\scriptstyle (b,k)\in H_l\atop
\scriptstyle t_bm > t_a k}
(1-v_{k}\sab)^{(\alpha_a|\alpha_b)(t_bm-t_ak)}.
\end{align}
The Jacobian $\partial z / \partial v $ is 1
at $v = 0$, so that
the map $z(v)$
is biholomorphic around $v=z=0$.
Let $v=v(z)$ the inverse map around $z=0$.

\begin{lem}[\cite{K2,HKOTY}]\label{lem:hkoty2}
Let $\beta\am\ ((a,m)\in H_l)$ be arbitrary
complex numbers.
We have the following power series
expansion at $z=0$ which converges for
$|z\am|<1$.
\begin{align}\label{eq:psi3}
\begin{split}
&\prod_{(a,m)\in H_l}
\left(1-v\am(z)\right)^{-\beta\am - 1}\\
&\qquad
 = \sum_{N\in \cN_l} \prod_{(a,m)\in H_l}
\binom{\beta\am +
c\am
 + N\am}{N\am}
\left(z_{m}\sa\right)^{N\am},
\end{split}
\end{align}
\begin{align}\label{eq:cma2}
c\am = \, & \sum_{\scriptstyle (b,k)\in H_l \atop
\scriptstyle t_b m < t_a k }(\alpha_a|\alpha_b)(t_a k-
t_b m)N\bk.
\end{align}
\end{lem}

A proof of Lemma \ref{lem:hkoty2} is given in Appendix 
\ref{sec:prhkoty}
 for reader's convenience.
\par
\subsection{Analytic formula}\label{subsec:analytic}

Let $w= (w\am)_{(a,m)\in H_l}$ be
 another complex multivariable.
We define
a biholomorphic map $w=w(v)$ around $v=w=0$ as
\begin{align}
w\am (v)& = v\am\prod_{(b,k)\in H_l}(1-v\bk
)^{-(\alpha_a|\alpha_b)\min(t_b m,t_a k)}\label{eq:wvrel}.
\end{align}
Combining it with (\ref{eq:zvrel}),
we have biholomorphic maps among
variables $v$, $z$, and $w$ around $v=z=w=0$.
Each map is denoted by $v=v(w)$, $z=z(w)$, {\it etc}.
By (\ref{eq:zvrel}) and (\ref{eq:wvrel}),
we have the relation
\begin{align}
z\am(w) &= w\am \prod_{(b,k)\in H_l}
(1-v\bk(w))^{(\alpha_a|\alpha_b)t_b m}.
\label{eq:wvrel2}
\end{align}

\begin{thm}\label{thm:RKgenerating}
The following equalities  (power series
expansions) hold around $w=0$: 
\begin{align}
\cK^\nu_l(w) &= \cK^0_l(w)
 \prod_{(a,m)\in H_l}(1-v\am(w))^{-\gamma\am}
,\label{eq:K1}\\
\cK^0_l(w) &=
\det_{H_l}
\left( 
\frac{ w\bk}{ v\am}
\frac{\partial v\am}{\partial w\bk}(w)
\right)
\prod_{(a,m)\in H_l}(1-v\am(w))^{-1}
\label{eq:K2}\\
&=\left(
\det_{H_l}
\left(\delta_{ab}\delta_{mk}+D_{am,bk}v\am(w)\right)
\right)^{-1},
\label{eq:R2}\\
\cR^\nu_l(w) &=  \prod_{(a,m)\in H_l}(1-v\am(w))^{-\gamma\am},
\label{eq:R1}
\end{align}
where $\gamma\am$ and  $D_{am,bk}$ are defined in
(\ref{eq:gammadef2}) and (\ref{eq:Dmat}).
\end{thm}
\begin{proof}
(\ref{eq:K1}) and (\ref{eq:K2}). Let $dz$ abbreviate
$\bigwedge_{(a,m)\in H_l}d z\am $.
By Lemma \ref{lem:hkoty2}, we have
\begin{align}
\label{eq:basic}
\begin{split}
&\quad \prod_{(a,m)\in H_l}
\binom{\beta\am +c\am + N\am}{N\am}\\
&= \Res_{z = 0}\left(
\prod_{(a,m)\in H_l}(1-v\am(z))^{-\beta\am-1}(z\am)^{-N\am-1} \right)
{dz}
\end{split}
\end{align}
In (\ref{eq:basic}), we set
 $\beta\am =P\am(\nu,N)-c\am=
 \gamma\am - \sum_{(b,k)\in H_l}
(\alpha_a|\alpha_b)t_akN\bk$.
After the substitution of  (\ref{eq:wvrel2}) and
the change of the integration variable, we obtain
\begin{align}
\label{eq:basic2}
\begin{split}
&\quad\prod_{(a,m) \in H_l}\binom{P\am + N\am}{N\am}\\
&= \Res_{w=0}
\Biggl\{
\Biggl(
\prod_{(a,m)\in H_l}(1-v\am(w))^{-\gamma\am-1}
(w\am)^{-N\am-1}
\Biggr)\\
&\qquad\qquad
\times 
\det_{H_l}\left(
\frac{w\bk}{z\am}
\frac{\partial z\am}{\partial w\bk}(w)\right)
\Biggr\}
{dw}.
\end{split}
\end{align}
Also, by (\ref{eq:zvrel}), we have
\begin{align}
\label{eq:detwz1}
\det_{H_l}
\left(
\frac{z\bk}{v\am}
\frac{\partial v\am}{\partial z\bk}
\right)=1,
\end{align}
The equalities (\ref{eq:K1}) and (\ref{eq:K2})
follows from 
(\ref{eq:basic2}), (\ref{eq:detwz1}),
 and the fact $\gamma=0$ if $\nu=0$.

\par
(\ref{eq:R2}).
By using  (\ref{eq:wvrel}),
the RHS of (\ref{eq:K2}) is easily calculated
as (\ref{eq:R2}).

\par
(\ref{eq:R1}).
In  (\ref{eq:basic2}), replace $N\am$
and $P\am$ in the both hand sides
 by $N[J]\am$ and $P[J]\am$.
Accordingly, $\gamma\am$ in the RHS
in  (\ref{eq:basic2})
should be  also replaced
by $\gamma\am - \sum_{(b,k)\in J}
(\alpha_a|\alpha_b)\min(t_b m,t_a k)$ (cf.\
(\ref{eq:peq}) and (\ref{eq:bsum})).
Then, using (\ref{eq:wvrel}), we obtain
\begin{align}
\label{eq:basic3}
\begin{split}
&\quad\prod_{(a,m) \in H_l}\binom{P[J]\am + N[J]\am}{N[J]\am}\\
&
=
\Res_{w=0}\left\{
\left(
\prod_{(a,m)\in H_l}
(1-v\am(w))^{-\gamma\am}
(w\am)^{-N\am-1}\right)\right.\\
&\qquad\qquad \times
\left.\left(\prod_{(a,m) \in J}v\am(w)\right)
K^0_l(w) \right\}
{dw}.
\end{split}
\end{align}
The equality (\ref{eq:R1})
follows from Lemma \ref{lem:ralt2},
(\ref{eq:R2}),
(\ref{eq:basic3}),
and the identity
\begin{align*}
\sum_{J \subset H_l}
\left(D_J
\prod_{(a,m) \in J}v\am(w)\right)=
\det_{(a,m), (b,k)\in H_l}
\left(\delta_{ab}\delta_{mk}+D_{am,bk}v\am(w)\right).
\end{align*}
\end{proof}
\begin{cor}\label{cor:Rlfact}
\begin{align}
\label{eq:rkkrel}
\cR^\nu_l(w) &= {\cK^\nu_l(w)}/{\cK^0_l(w)},\\
\cR^\nu_l(w)\cR^{\nu'}_l(w) &= \cR^{\nu+\nu'}_l(w).
\label{eq:factorization}
\end{align}
\end{cor}
\begin{proof}
(\ref{eq:rkkrel}) is immediately obtained from
(\ref{eq:R1}) and (\ref{eq:K1}).
(\ref{eq:factorization}) follows from
(\ref{eq:R1}) and the fact
$\gamma(\nu+\nu')\am = \gamma(\nu)\am + \gamma(\nu')\am$.
\end{proof}

\begin{rem}\label{rem:main}
Let $\mathbb{C}_l[[w]]$ be the 
(formal) power series ring of $(w\am)_{(a,m)\in H_l}$.
There are 
natural surjections $\varphi_{lk}:
\mathbb{C}_k[[w]]\rightarrow \mathbb{C}_l[[w]]$ 
($l \leq k$) with
$\varphi_{lk}(w\am)=0$ for $(a,m)\in H_{k}\setminus H_l$.
Since $\varphi_{lk}(\cR^\nu_k(w))=\cR^\nu_l(w)$,
$(\cR^\nu_l(w))_{l=1}^\infty$ defines
an element in 
the projective limit
$\mathbb{C}[[w]]:=\varprojlim\mathbb{C}_l[[w]]$,
which is identified with 
$\cR^\nu(w)$  in (\ref{eq:rwdef})
(and so is $\cK^\nu(w)$).
Then, the equalities (\ref{eq:R1}), (\ref{eq:K1}),
and (\ref{eq:factorization}) can be rephrased
as follows:
\begin{gather}\label{eq:RKKinf}
\cR^\nu(w) = {\cK^\nu(w)}/{\cK^0(w)}
= \prod_{(a,m)\in H}(1-v\am(w))^{-\gamma\am},\\
\cR^\nu(w)\cR^{\nu'}(w) = \cR^{\nu+\nu'}(w),
\end{gather}
where the RHS of (\ref{eq:RKKinf}) denotes
$(f_l(w))\in \bC[[w]]$ with $f_l(w)$ being
the power series expansion of the
the RHS of (\ref{eq:R1}) around $w=0$.
\par
\end{rem}

\section{Formal completeness of Bethe vectors}

Now we are ready to present
the main results briefly explained
in Section \ref{subsec:xxz}.
Our goal here is to express the coefficient
$r^\nu_\lambda$ in (\ref{eq:rres1})
as a sum of the numbers  $R(\nu,N)$.

\subsection{Specialization of generating series}
\label{subsec:special}

We remind that $y_a=e^{-\alpha_a}$ and
$x_a=e^{\Lambda_a}$
are the formal
exponents of simple roots and fundamental weights of
$X_n$.
We also regard $y=(y_a)$ and $x=(x_a)$
as complex multivariables related by the map
$y=y(x)$, 
$y_a=\prod_{b=1}^n x_b^{-(\alpha_a|\alpha_b)t_b}$.
We do the specialization  $w\am(y)=y^m_a$
of the variable
of the series $\cR^\nu_l(w)$,
\begin{align}
\cR^\nu_l(w(y))
&=
\sum_{N\in \cN_l} R(\nu,N) \prod_{a=1}^n
y_a^{ \sum_{m=0}^\infty m N\am}.
\end{align}
The limit 
\begin{align}
\label{eq:tr5}
\tR^\nu (y):=&\,
\lim_{l \to \infty}
\cR^\nu_l(w(y))
=
\sum_{N\in \cN} R(\nu,N) \prod_{a=1}^n
y_a^{ \sum_{m=0}^\infty m N\am}
\end{align}
exists in $\mathbb{C}[[y]]$,
because
$\cR^\nu_l(w(y))\equiv
\tR^\nu(y)$ mod $I_l$,
where $I_l$ is the ideal of $\bC[[y]]$
generated by
$y_1^{t_1l+1}$, \dots,
$y_n^{t_nl+1}$.
In the same way, we define the limit
\begin{align}
\label{eq:tk5}
\tK^\nu (y):=
&\,
\lim_{l \to \infty}
\cK^\nu_l(w(y))
=\sum_{N\in \cN} K(\nu,N) \prod_{a=1}^n
y_a^{ \sum_{m=0}^\infty m N\am}.
\end{align}
For each $(a,m)\in H$,
let $\delta\am=(\nu\bk)_{(b,k)\in H}$, $\nu\bk=\delta_{ab}
\delta_{mk}$ and
\begin{align}\label{eq:ram1}
\tR\am(y):=\tR^{\delta\am}(y),\quad
\tK\am(y):=\tK^{\delta\am}(y).
\end{align}
\par
It immediately follows from Corollary \ref{cor:Rlfact} that
\begin{prop}\label{prop:RK1}
\begin{gather}
\label{eq:RKK2}
\tR^\nu(y) = {\tK^\nu}(y)/{\tK^0}(y),
\quad \tR\am(y)=\tK\am(y)/\tK^0(y),\\
\label{eq:fac2}
\tR^\nu (y)\tR^{\nu'}(y) = \tR^{\nu+\nu'}(y),\\
\label{eq:fac3}
\tR^\nu(y)=\prod_{(a,m)\in H} (\tR\am(y))^{\nu\am},\\
\label{eq:KK1}
\tK^\nu (y)/ \tK^0(y)=\prod_{(a,m)\in H} 
(\tK\am (y)/ \tK^0(y))^{\nu\am}.
\end{gather}
\end{prop}
We also introduce the corresponding Laurent series of $x$ as follows:
\begin{alignat}{2}
R\am(x) &= x_a^m \tR\am(y(x)),\quad &
R^\nu(x)&
=\biggl(\prod_{(a,m)\in H} x_a^{m\nu\am}\biggr)
\tR^\nu(y(x)),\\
\label{eq:kam}
K\am(x) &= x_a^m \tK\am(y(x)),\quad&
K^\nu(x)&=
\biggl(\prod_{(a,m)\in H} x_a^{m\nu\am}\biggr)
\tK^\nu(y(x)).
\end{alignat}
\begin{prop}\label{prop:RK2}
The equalities in (\ref{eq:RKK2})--(\ref{eq:KK1}) with
$\tR^\nu$, $\tR\am$, $\tK^\nu$, $\tK\am$
being replaced with
$R^\nu$, $R\am$, $K^\nu$, $K\am$
 also hold.
\end{prop}

\subsection{Main Results}

\begin{thm}\label{thm:main3}
There exists
a unique
 family $(\tQ\am)_{(m,a)\in H}$ of invertible
 power series of $y$
which
satisfies (\~Q-I) and (\~Q-II) in
Definition \ref{defn:unique}. 
In fact, 
\begin{align}
\tQ\am(y)=\tR\am(y)=\tK\am(y)/\tK^0(y).
\end{align}
\end{thm}

\begin{rem}
The existence of $(\tQ\am)_{(m,a)\in H}$ also follows
from Theorem \ref{thm:hkoty1}
for $X_n$ of classical type.
A  weak version of the uniqueness was shown in 
\cite[Theorem 8.1]{HKOTY},
where the $\cW$ invariance of $Q\am$ was further assumed.
\end{rem}

A proof of Theorem \ref{thm:main3}
is given in Section \ref{subsec:proofmain}.
We  present its consequences first.
Let $Q^\nu(x)$ be the  Laurent series of $x$
defined in (\ref{eq:Qnu}).
It follows from Propositions \ref{prop:RK1}, \ref{prop:RK2}
 and (\ref{eq:Qnu})
that
\begin{cor}
\begin{align}
\tQ^\nu(y) &= \tR^\nu(y) = \tK^\nu(y)/\tK^0(y),\\
\label{eq:qrk1}
Q^\nu(x) &= R^\nu(x) = K^\nu(x)/K^0(x).
\end{align}
\end{cor}

Expanding the both hand sides
of the first equality in (\ref{eq:qrk1}),
we have
\begin{align}\label{eq:rwt2}
r_\lambda^\nu
=\sum_{N\in \cN^\nu_\lambda}R(\nu,N),
\end{align}
where $r_\lambda^\nu$ 
is defined in (\ref{eq:Qexp5}), and
\begin{align}
\label{eq:nset1}
\cN^\nu_\lambda &=
\{\,N\in \cN\mid
\sum_{(a,m)\in H}m\nu\am\Lambda_a-
\sum_{(a,m)\in H}m N\am\alpha_a
=\lambda\,\}.
\end{align}
Therefore,

\par
\begin{cor}[Formal completeness of \XXZ-type
 Bethe vectors]\label{cor:completeness3}
If Conjecture
\ref{conj:qsys} is correct,
then
\begin{align}
 \mathrm{ch}\,W^\nu (x)&=R^\nu (x), \\
\label{eq:rwt1}
 \dim W^\nu_\lambda &=
\sum_{N\in \cN_\lambda^\nu}R(\nu,N),
\end{align}
where $\dim W^\nu_\lambda$ denotes the weight multiplicity
 in $W^\nu$
 at weight $\lambda$.
\end{cor}

\begin{rem}
We know that the number $R(\nu,N)$ 
correctly counts the Bethe vectors
only for special string patterns $N$.
One  naive explanation of the equality
(\ref{eq:rwt1})  is as follows:
$R(\nu,N)$ correctly counts the Bethe vectors,
therefore, the weight multiplicity when
$\lambda$ are relatively close to the highest weight of
$W^\nu$.  Then, the factorization property
(\ref{eq:fac3}) imposes such a strong constraint that
the equality (\ref{eq:rwt1}) has to hold for the entire
 region of $\lambda$.
\end{rem}

\par
With Theorem \ref{thm:hkoty1}, we have
\begin{cor} Let $X_n$ be of classical type,
and let $\chi\am$ be the $X_n$-character
given in the RHSs of
(\ref{eq:a1})--(\ref{eq:d1}).
Then
\begin{align}
\chi\am(x)=R\am(x)={K\am}(x)/{K^0}(x).
\end{align}
In particular, 
$R^\nu(x)=K^\nu(x)/K^0(x)$ is a $\cW$-invariant
Laurent polynomial, and
$\tR^\nu(y)=\tK^\nu (y)/ \tK^0(y)$ is a polynomial.

\end{cor}

\begin{rem}
It remains an open problem to show 
 the $\cW$-invariance of $Q^\nu(x)$ and the
polynomial property of $\tQ^\nu(y)$
for $X_n$ of exceptional type
without assuming Conjecture \ref{conj:qsys}.
\end{rem}

The formal completeness of the \XXX-type Bethe vectors
has been worked out 
by \cite{K1, K2, KR, HKOTY}.
We reformulate their result in 
our context as follows (See Appendix \ref{sec:ksec} for  a proof):

\begin{prop}\label{prop:weyl1}
If $Q_1^{(1)}(x)$,\dots, $Q_1^{(n)}(x)$  are $\cW$-invariant,
then 
\begin{align*}
K^0(x)
=\prod_{\alpha \in \Delta_+} (1-e^{-\alpha}),
\end{align*}
where $\Delta_+$ is the set of all the positive
roots of $X_n$.
\end{prop}
If Conjecture \ref{conj:qsys} is correct,
then $Q_1^{(1)}(x)$,\dots, $Q_1^{(n)}(x)$ are $\cW$-invariant.
Then, by Proposition \ref{prop:weyl1}, we have
\begin{align}\label{eq:kwt2}
k_\lambda^\nu
&=\sum_{N\in \cN_\lambda^\nu}K(\nu,N),
\end{align}
where $k_\lambda^\nu$ is defined in (\ref{eq:Qexp4}).
Therefore,
\begin{cor}[Formal completeness of \XXX-type
 Bethe vectors]\label{cor:completeness4}
If Conjecture \ref{conj:qsys} is correct,
then
\begin{align}
  \mathrm{ch}\,W^\nu (x)&=
K^\nu(x)/ \prod_{\alpha \in \Delta_+}(1-e^{-\alpha}),\\
 [W^\nu: V_\lambda]
&=
\sum_{N\in \cN_\lambda^\nu}
 K(\nu,N),
\end{align}
where $\lambda$ is a dominant  $X_n$-weight,
and 
$[W^\nu : V_\lambda]$ denotes the multiplicity 
of the $\Uqcl$-irreducible components $V_\lambda$
with highest weight $\lambda$ in $W^\nu$.
\end{cor}

\subsection{Proof of Theorem \ref{thm:main3}}
\label{subsec:proofmain}
\subsubsection{Existence}
Because $R(\nu,N=0)=1$, $\tR\am(y)$ is an invertible
power series.
We show that $(\tR\am)$ satisfies (\~Q-I) and (\~Q-II)
with $\tQ\am$ being replaced with $\tR\am$.
\par
(\~Q-I).
{}From the factorization property (\ref{eq:fac3}),
it is enough to prove 
\begin{lem}\label{lem:sum}
The following relation holds:
\begin{equation}\label{eq:qrel1}
\tR^{\lambda}(y) = 
\tR^{\mu} (y)+ 
y_a^m\tR^{\nu}(y),
\end{equation}
where $\lambda = (\lambda\bk)$,
$ \mu = (\mu\bk)$, $\nu = (\nu\bk)$ with
\begin{align*}
\lambda\bk &= 2\delta_{ab}\delta_{mk},\quad
\mu\bk = \delta_{ab}(\delta_{m+1,k}+
\delta_{m-1,k}),\\
\nu\bk&=2\delta_{ab}\delta_{mk}-(\alpha_a|\alpha_b)
B_{am,bk}.
\end{align*}
\end{lem}
\begin{proof}
It is enough to show that
\begin{equation}\label{eq:req}
R(\lambda,N) = R(\mu,N) + R(\nu,N')
\end{equation}
for $N=(N\bk)\in \tilde{\cN}$
 and $N'=(N_k^{\prime(b)}) \in \tilde{\cN}$ which are related
as $N_k^{\prime (b)}=N\bk - \delta_{ab}\delta_{mk}$
($\tilde{\cN}$ is defined in (\ref{eq:tcN})).
By Remark \ref{rem:Next}, for any $N\in \tilde{\cN}$ it holds that
\begin{align}
\label{eq:rf1}
 R(\lambda,N) 
= \sum_{J \subset H}
\left\{ D_J
\prod_{(b,k) \in H}\binom{P\bk(\lambda[J],N[J])
 + N[J]\bk}{N[J]\bk}
\right\}.
\end{align}
By (\ref{eq:bformula})
and (\ref{eq:bsum}), it is easy to show
\begin{align}
P\bk(\lambda[J],N[J]) = P\bk(\mu[J],N[J]) +
\delta_{ab} \delta_{mk}\label{eq:ps}
= P\bk(\nu[J],N'[J]).
\end{align}
With (\ref{eq:ps}), we have
\begin{align}
\label{eq:pe1}
\begin{split}
\binom{P\am(\lambda[J],N[J])  + N[J]\am}{N[J]\am}
&=
\binom{P\am(\mu[J],N[J]) + N[J]\am}{N[J]\am}\\
&\qquad
+\binom{P\am(\nu[J],N'[J]) + N'[J]\am}{N'[J]\am},
\end{split}
\end{align}
while for $(b,k)\neq (a,m)$,
\begin{align}
\label{eq:pe2}
\begin{split}
\binom{P[J]\bk + N[J]\bk}{N[J]\bk}
&=
\binom{P\bk(\mu[J],N[J]) + N[J]\bk}{N[J]\bk}\\
&
=\binom{P\bk(\nu[J],N'[J]) + N'[J]\bk}{N'[J]\bk}.
\end{split}
\end{align}
The equality (\ref{eq:req}) follows from
(\ref{eq:rf1}), (\ref{eq:pe1}), and (\ref{eq:pe2}).
\end{proof}

\par
(\~Q-II).
We show the limit
$\lim_{m\to \infty} \tR\am(y)$ exists in $\mathbb{C}
[[y]]$.
Let $\delta\am$ be the one in (\ref{eq:ram1}).
Then, $P\bk(\delta\am,N)
 = P\bk(\delta_{m+1}\sa,N) - 
\delta_{ab}\theta(k \ge
m+1)$ holds from (\ref{eq:pdef2}).
In the series $\tR\am(y)$,
 those $N = (N\bk)$ containing 
$N_k\sa >0$ with  $k \ge m+1$  make contribution
to the power $y_a^d$ only for $d>m$ (see (\ref{eq:tr5})).
It follows that
$\tR\am (y)\equiv \tR_{m+1}\sa(y)$ mod
 $y_a^{m+1}\mathbb{C}[[y]]$.
Then, we have
\begin{displaymath}
\tR\am (y)\equiv \tR_{m+1}\sa (y)\equiv
\tR_{m+2}\sa (y)\equiv \cdots
\quad \text{mod}\ y_a^{m+1} \mathbb{C}[[y]],
\end{displaymath}
which means $\lim_{m\to \infty} \tR\am(y)$ exists.

\subsubsection{Uniqueness}
\label{subsubsec:unique}

Let $(\tQ\am(y))_{(a,m)\in H}$ be 
 a family of invertible power series
of $y$ which
 satisfies (\~Q-I) and
(\~Q-II) in
Definition \ref{defn:unique}.
By Proposition \ref{prop:rec1} and
(\~Q-II), the constant term
of $\tQ\am(y)$ is $1$.
We define a family of power series
$(v\am(y))_{(a,m)\in H_l}$
with constant term zero
by 
\begin{equation}\label{eq:vam}
v\am(y) = 1-\frac{\tQ_{m-1}\sa(y) \tQ_{m+1}\sa(y)}
{(\tQ\am(y))^2},
\quad \tQ_0\sa(y)=1.
\end{equation}
We further define a  family of
power series
$(w\am(y))_{(a,m)\in H_l}$
with constant term zero
by the composition 
$w\am(v(y))$ of the series $v\am(y)$
and the power series expansion of
the holomorphic map
$w=w(v)$ in (\ref{eq:wvrel})
at $v=0$.

\begin{lem}\label{lem:spe1}
The following equalities of power series of $y$ hold:
\begin{align}
\prod_{(a,m)\in H_l}(1-v\am(y))^{-\gamma\am}
&=
\label{eq:vs}
\left(
\prod_{(a,m)\in H_l}
(\tQ\am(y))^{\nu\am}
\right)
\left(
\prod_{a=1}^n\frac{(\tQ_{t_a l}\sa(y))^{\gamma_{t_al+1}\sa}}
{(\tQ_{t_a l+1}\sa(y))^{\gamma_{t_al}\sa}}
\right),\\
\label{eq:ws}
w\am(y) &= y_a^m \prod_{b=1}^n
\left(\frac{\tQ_{t_b l}\sab(y)}{\tQ_{t_b l+1}\sab(y)}
\right)^{(\alpha_a|\alpha_b)t_bm}.
\end{align}
\end{lem}

\begin{proof}
(\ref{eq:vs}). By rearranging the product indices, 
the LHS in (\ref{eq:vs}) becomes
\begin{equation*}
\prod_{a=1}^n
\left(
\frac{(\tQ_{t_a l}\sa(y))^{\gamma_{t_al+1}\sa}}
{(\tQ_{t_a l+1}\sa(y))^{\gamma_{t_al}\sa}}
\prod_{m=1}^{t_al}(\tQ\am(y))^{%
2\gamma\am - \gamma_{m+1}\sa-\gamma_{m-1}\sa}
\right).
\end{equation*}
By (\ref{eq:minrel1}), we have
$2\gamma\am - \gamma_{m+1}\sa-\gamma_{m-1}\sa=\nu\am$.
\par
(\ref{eq:ws}). Using the same trick as above and 
the definition of $B_{am,bk}$ in (\ref{eq:bfunc}), we have
\begin{align*}
&\prod_{(b,k)\in H_l}(1-v\bk(y))^{%
-(\alpha_a|\alpha_b)\min(t_bm,t_ak)}\\
=&\prod_{b=1}^n
\left(
\left(\frac{\tQ_{t_b l}\sab(y)}{\tQ_{t_b l+1}\sab(y)}
\right)^{(\alpha_a|\alpha_b)t_bm}
\prod_{k=1}^{t_bl}(\tQ\bk(y))^{%
(\alpha_a|\alpha_b)B_{am,bk}}
\right).
\end{align*}
On the other hand, using (\~Q-I) and 
(\ref{eq:b0cond1}), we have
\begin{align*}
v\am (y)=
 y_a^m \prod_{(b,k)\in H_l}
(\tQ\bk(y))^{-(\alpha_a|\alpha_b)B_{am,bk}}.
\end{align*}
for $(a,m)\in H_l$.
(\ref{eq:ws}) is obtained by multiplying the above two equalities.
\end{proof}

\begin{lem}\label{lem:qrel1}
The following equality of power series of $y$ holds:
\begin{align}
\label{eq:vs3}
\begin{split}
&\quad 
\left(
\prod_{(a,m)\in H_l}
(\tQ\am(y))^{\nu\am}
\right)
\left(
\prod_{a=1}^n\frac{(\tQ_{t_a l}\sa(y))^{\gamma_{t_al+1}\sa}}
{(\tQ_{t_a l+1}\sa(y))^{\gamma_{t_al}\sa}}
\right)\\
&=
\sum_{N\in \cN_l}
R(\nu,N)
\prod_{(a,m)\in H_l}
\left(
 y_a^{mN\am}
\prod_{b=1}^n
\left(
\frac{(\tQ_{t_b l}\sab(y))}
{(\tQ_{t_b l+1}\sab(y))}
\right)^{(\alpha_a|\alpha_b)t_bm N\am}
\right).
\end{split}
\end{align}
\end{lem}
\begin{proof}
Let us regard
(\ref{eq:R1})
as an equality of power series of $w$.
Then, by substituting the
series $w\am(y)$ for the variable $w\am$
in (\ref{eq:R1}) and
using  (\ref{eq:vs}) and (\ref{eq:ws}),
we obtain (\ref{eq:vs3}).
\end{proof}

\par
Thanks to the convergence property (\~Q-II),
the limit $l\rightarrow \infty$ of (\ref{eq:vs3})
gives the equality
$ \prod_{(a,m)\in H} (\tQ\am(y))^{\nu\am}=
\tR^\nu(y)$.
In particular, by setting $\nu=\delta\am$, we obtain
$\tQ\am(y)= \tR\am(y)$.
This completes the proof of
the uniqueness property of $(\tQ\am)$,
thereby finishes the proof of Theorem \ref{thm:main3}.


\appendix

\section{Some properties of the function $B_{am,bk}$}
\label{sec:Bfunc}

The function $B_{am,bk}$, 
\begin{equation}
B_{am,bk}:=
2\min(t_bm,t_ak)-\min(t_bm,t_a(k+1))
-\min(t_bm,t_a(k-1))
\end{equation}
appears in several places such as
(\ref{eq:qsys6}),
(\ref{eq:qsys1}),
(\ref{eq:nuJ}), (\ref{eq:qrel1}),
(\ref{eq:ws}), {\em etc},
and plays a key role.
Below we list  the properties we use.
Since it always appears in the combination
$(\alpha_a | \alpha_b) B_{am,bk}$,
we are interested only in the situations
$(t_a,t_b)=(1,1)$, $(2,2)$, $(3,3)$,
 $(1,2)$, $(1,3)$, $(2,1)$, $(3,1)$
here.

Let us first observe that the infinite-size
matrix $B=(B_{am,bk})_{(a,m),(b,k)\in H}$
is expressed as a product $B=B'D$,
where $B'=(B_{am,bk})$, $D=(D_{am,bk})$
with
\begin{align}
B'_{am,bk}&=\min(t_bm,t_ak),\\
D_{am,bk}&=\delta_{ab}(2\delta_{mk}
-\delta_{m,k-1}-\delta_{m,k+1}).
\end{align}
Since the relation
\begin{align}\label{eq:minrel1}
2\min(m,k)-\min(m,k+1)-\min(m,k-1)=\delta_{mk}
\end{align}
holds, the inverse matrix $D^{-1}$ of   $D$ is given by
\begin{align}
(D^{-1})_{am,bk}=\delta_{ab}\min(m,k).
\end{align}

\begin{prop}
\label{prop:bfunc1}
\par
(i) For each $(a,m)$, there are only finitely-many
$(b,k)$'s such that $B_{am,bk}\neq 0$.
Explicitly, 
\begin{equation}\label{eq:bformula}
B_{am,bk}=
\begin{cases}
2\delta_{m,2k}+\delta_{m,2k+1}+\delta_{m,2k-1}
& (t_a,t_b)=(2,1)\\
3\delta_{m,3k}+2\delta_{m,3k+1}+2\delta_{m,3k-1}&
 (t_a,t_b)=(3,1)\\
\qquad\qquad
+\delta_{m,3k+2}+\delta_{m,3k-2}&\\
t_a\delta_{t_bm,t_a k}
& \text{otherwise}.
\end{cases}
\end{equation}
(ii)
Let $H_l$
 be the subset of $H$ 
defined in (\ref{eq:hl1}).
Then,
\begin{align}
\label{eq:b0cond1}
B_{am,bk}=0\quad
\text{for $(a,m)\in H_l$, $(b,k)\notin H_l$}.
\end{align}
(iii) The following relations hold:
\begin{align}
\label{eq:bsum}
\sum_{j=1}^\infty
B_{am,bj}\min(j,k)
&=\min(t_bm,t_ak),\\
\label{eq:bfunc5}
\sum_{k=1}^\infty 
B_{am,bk}k &=
t_b m,\\
\label{eq:bfunc4}
\sum_{(b,k)\in H}
(\alpha_a|\alpha_b)B_{am,bk}
k \Lambda_b &= m\alpha_a.
\end{align}
\end{prop}

\begin{proof}
(i) This is shown
by the case check.
(ii) This can be easily checked by
 (\ref{eq:bformula}).
(iii) (\ref{eq:bsum})
is equivalent to the matrix relation
 $(B'D)D^{-1}=B'$.
We have only to care  that
the matrix product in the LHS is well-defined.
This is guaranteed  by (i).
The LHS of (\ref{eq:bfunc5}) can be calculated
 in a similar way as follows:
\begin{align*}
\sum_{k=1}^L
\biggl( \sum_{j=1}^\infty
B'_{am,bj}D_{bj,bk}\biggr)k
=
(L+1)B'_{am,bL}-LB'_{am,bL+1}
=t_bm,
\end{align*}
where  $L$ is a sufficiently large number.
(\ref{eq:bfunc4}) immediately follows from (\ref{eq:bfunc5})
and the relation $\alpha_a=\sum_{b=1}^n
(\alpha_a|\alpha_b)t_b \Lambda_b$.
\end{proof}

As an application of  (\ref{eq:bformula}),
we show that
\begin{prop}
\label{prop:rec1}
The relation (\~Q-I) 
recursively
determines all the other power
series $\tQ\am(y)$ ($m\geq 2$) from given invertible
power series $\tQ_1^{(1)}(y)$,
\dots, $\tQ_1^{(n)}(y)$ as an initial condition;
furthermore,
so determined power series $\tQ\am(y)$ is invertible,
and its constant term $c\am$ is the $m$-th power
of the constant term of\/ $\tQ_1\sa(y)$.
\end{prop}
\begin{proof}
We introduce another subset $\cH_l$ ($l\geq 1$) of the
index set $H$ as
\begin{align}
\label{eq:hset2}
\cH_l=\{\, (a,m)\in H \mid t(m-1)\leq t_a(l-1)\, \},
\quad
t=\max_{1\leq a \leq n} t_a.
\end{align}
Then, $\cH_1=\{\, (a,1)\mid 1\leq a\leq n\, \}
\subset \cH_2\subset \cdots$ and $\varinjlim \cH_l = H$.
By (\~Q-I), 
\begin{equation}\label{eq:qsysalt2}
\tQ_{m+1}\sa (y)
=
\frac{(\tQ\am(y))^2}{\tQ_{m-1}\sa(y)}
\Biggl(
1-y_a^m
\prod_{(b,k)\in H}
(\tQ\bk(y))^{-(\alpha_a|\alpha_b)B_{am,bk}}
\Biggr).
\end{equation}
For a given $(a,m+1)\in H$,
let $l$ be a unique positive integer such that
$(a,m+1)\in \cH_{l+1}\setminus \cH_l$.
Then, with (\ref{eq:bformula}),
it is easy to check that
$B_{am,bk}=0$ for $(b,k)\notin \cH_l$.
The claim now follows from (\ref{eq:qsysalt2})
by induction on $l$.
\end{proof}

\section{Proof of Lemma \ref{lem:hkoty2}}\label{sec:prhkoty}

The following proof
of Lemma \ref{lem:hkoty2}
 is essentially quoted from \cite[Proposition 8.3]{HKOTY}.

Let $t:=\max\{\, t_1,\dots,t_n\, \}$ and
\begin{eqnarray}
H_l[i] & = & \{\, (a,m) \in H_l\mid  t m\ge t_a i \,\}.
\label{eq:Hdef1}
\end{eqnarray}
Then $\emptyset=H_l[tl+1]\subset
H_l[tl]\subset \dots \subset  H_l[1] = H_l$.
For each $1\leq i \leq  tl$,
let
$z_{i}=z_{i}(v_i)$,
$v_i=(v_{m,i}\sa)_{(a,m)\in H_l[i]}$,
 $z_{i}=(z_{m,i}\sa)_{(a,m)\in H_l[i]}$
be the biholomorphic map
 around $v_i=z_{i}=0$ 
defined by
\begin{equation}\label{eq:zsystem}
z_{m,i}\sa (v_i)= v_{m,i}\sa
\prod_{(b,k)\in H_l[i]\setminus H_l[tm/t_a]}
\left(1-v_{k,i}\sab
 \right)^{(\alpha_a|\alpha_b)(t_b m-t_a k)},
\end{equation}
and $v_i=v_i(z_{i})$ be its inverse.
($z_{m,i}\sa$ here corresponds to $z_{m,i-1}\sa$ in
\cite{HKOTY}.)
Let $v_{i+1}=v_{i+1}(v_i)$ be the holomorphic map defined by
$v_{m,i+1}\sa(v_i)=v_{m,i}\sa$ (for $(a,m)\in H_l[i+1]$),
and $z_{i+1}=z_{i+1}(z_i)$ be the
composition $z_{i+1}(v_{i+1}(v_i(z_i)))$.
Namely,
\begin{equation}\label{eq:zsystem2}
z_{m,i+1}\sa (z_{i})= z_{m,i}\sa
\prod_{(b,k)\in H_l[i]\setminus H_l[i+1]}
\left(1-v_{k,i}\sab(z_{i})
 \right)^{-(\alpha_a|\alpha_b)(t_b m-t_a k)}.
\end{equation}
The relation of these variables and maps
 are summarized by the following
diagram:
\begin{align}
\begin{matrix}
v_{i+1} & \leftarrow & v_i\\
\updownarrow & & \updownarrow\\
z_{i+1} & \leftarrow & z_i \\
\end{matrix}
\, .
\end{align}

The condition $(b,k)\in H_l[i]\setminus H_l[tm/t_a]$
is equivalent to $(b,k)\in H_l[i]$ and $t_bm>t_ak$.
Thus, if we set $z_1=z$ and $v_1=v$,
the map $z_1(v_1)$ coincides with $z(v)$ in (\ref{eq:zvrel}).
Lemma \ref{lem:hkoty2} is a special case $i=1$ of the following
proposition.
\begin{prop}[\cite{K2,HKOTY}]\label{lem:hkoty}
For any integer $1 \le i \le tl$,
and any complex numbers $\beta\am$ $((a,m)\in H_l[i])$,
we have the  following power series expansion at $z=0$
which converges  for $|z_{m,i}\sa|<1$:
\begin{align}\label{eq:psi2}
\begin{split}
&\prod_{(a,m)\in H_l[i]}
\left(1-v_{m,i}\sa(z_{i})\right)^{-\beta\am - 1}\\
&\qquad  = \sum_{ N\in \cN_l[i]} \prod_{(a,m)\in H_l[i]}
\binom{\beta\am +
c\am
 + N\am}{N\am}
\left(z_{m,i}\sa\right)^{N\am},
\end{split}
\end{align}
where
\begin{align}\label{eq:cma}
c\am =  & \sum_{ (b,k)\in H_l[tm/t_a+1] }
(\alpha_a|\alpha_b)(t_a k-
t_b m)N\bk,\\
\cN_l[i]=&
\{\, N=(N\am)_{(a,m)\in H} \mid
N\am \in \mathbb{Z}_{\ge 0},\
N\am = 0\ \text{for}\ (a,m)\notin H_l[i]\,\}.
\end{align}
\end{prop}

\begin{rem}

If $(a,m)\in H_l[i]$, then
$H_l[tm/t_a+1]\subset H_l[i]$. 
Also, the condition $(b,k)\in H_l[tm/t_a+1]$ is equivalent to
the condition $(b,k)\in H_l$, $t_bm < t_a k$.
Therefore, $c\am$ in (\ref{eq:cma}) is the same one as in
(\ref{eq:cma2}).
\end{rem}

\begin{proof}
We prove the proposition by induction on $i$ in the
descent order.
First, consider the case $i=tl$.
Suppose $(a,m)\in H_l[tl]=\{ (a,t_al)\}_{a=1}^n$.
Then,
$c\am=0$ due to $H_l[tm/t_a+1]
=H_l[tl+1]=\emptyset$,
and $v_{m,tl}\sa(z_{tl})=z_{m,tl}\sa$ due to
$H_l[tl]\setminus H_l[tm/t_a]=\emptyset$.
Therefore, the claim reduces to the well-known power series expansion
\begin{equation}\label{eq:exp}
(1-v)^{-\beta-1} = \sum_{N=0}^\infty \binom{\beta+N}{N}v^N,
\end{equation}
which converges for $|v|<1$.
Next, let us assume (\ref{eq:psi2}) holds for $i+1$.
Then, for $z_i$ such that
$|z_{m,i}\sa|<1$ and $|z_{m,i+1}\sa(z_i)|<1$,
the LHS of (\ref{eq:psi2})
is equal to
\begin{align}
\label{eq:psi1}
\begin{split}
& \left(\sum_{ N \in \cN_l[i+1]}
 \prod_{(a,m)\in H_l[i+1]}
\binom{\beta\am +
c\am
 + N\am}{N\am}
\left(z_{m,i+1}\sa(z_i)\right)^{N\am}\right)\\
& \qquad \times
\left(
\prod_{(a,m)\in 
H_l[i]\setminus H_i[i+1]}
 (1-v_{m,i}\sa(z_{i}))^{-\beta\am -1}
\right)
\end{split}
\end{align}
by the induction hypothesis.
Here, we used the fact
that, for $(a,m)\in H_l[i+1]$,
$v_{m,i}\sa(z_i)=v_{m,i+1}\sa(v_i(z_i))=
v_{m,i+1}\sa(z_{i+1}(z_i))$.
Substituting (\ref{eq:zsystem2}) for $z_{m,i+1}\sa (z_{i})$
in  (\ref{eq:psi1}), we have
\begin{align}\label{eq:psi5}
\begin{split}
& \sum_{ N \in \cN_l[i+1]}
\left(
 \prod_{(a,m)\in H_l[i+1]}
\binom{\beta\am +
c\am
 + N\am}{N\am}
\left(z_{m,i}\sa\right)^{N\am}\right)\\
& \qquad\qquad \times
\left(
\prod_{(a,m)\in 
H_l[i]\setminus H_i[i+1]}
 (1-v_{m,i}\sa(z_{i}))^{-\beta\am -\tilde{c}_m^{\sa}-1}
\right),
\end{split}
\end{align}
where
\begin{equation}
\tilde{c}\am =
\sum_{(b,k)\in H_l[i+1]} (\alpha_a|\alpha_b)
(t_a k - t_b m)N\bk.
\end{equation}
Suppose $(a,m)\in
H_l[i] \setminus H_l[i+1]$.
Then, ${tm}/{t_a}= i$ holds.
It follows that $\tilde{c}\am=c\am$ and
$v_{m,i}\sa(z_{i})=z\sa_{m,i}$.
Thus, applying (\ref{eq:exp}) to the second factor
of (\ref{eq:psi5}), we obtain (\ref{eq:psi2}).
It is easy to check that the RHS in (\ref{eq:psi2}) 
converges for $|z_{m,i}\sa|<1$.
\end{proof}

\section{Proof of Proposition \ref{prop:weyl1}}
\label{sec:ksec}

Following \cite{HKOTY},
we prove the proposition in two steps as
Propositions \ref{prop:k01} and \ref{prop:qx5}.
The proof of Proposition \ref{prop:qx5} is
taken from \cite{HKOTY},
while
the proof of Proposition \ref{prop:k01} here is new.

\par
{\em Step 1.}
We start from the formula (\ref{eq:K2}),
which is also written as (cf.\ (\ref{eq:detwz1}))
\begin{align}
\cK^0_l(w) &=
\det_{H_l}
\left( 
\frac{ w\bk}{ z\am}
\frac{\partial z\am}{\partial w\bk}(w)
\right)
\prod_{(a,m)\in H_l}(1-v\am(w))^{-1},
\label{eq:K3}
\end{align}
where $w\am$, $v\am$, and $z\am$ are related 
by (\ref{eq:zvrel}), (\ref{eq:wvrel}),
and (\ref{eq:wvrel2}).
Let $w\am(y)=y_a^m$ be the specialization in
Section \ref{subsec:special}.
We define new series of $y$
\begin{align}
f_a(y)=\prod_{m=1}^{t_al}(1-v\am(w(y))),
\quad a=1,\, \dots,\, n,
\end{align}
which depend also on $l$.

\begin{lem}
\label{lem:flim1}
In\/ $\bC[[u]]$,
\begin{align}
\lim_{l\to \infty}f_a(y)=
(\tQ\sa_1(y))^{-1}.
\end{align}
\end{lem}

\begin{proof}
By (\ref{eq:R1}),
$f_a(y)=\cR^\nu_l(w(y))$,
where $\nu=(\nu\bk)$ with
$\nu\bk=-1$ if $(b,k)=(a,1)$, $(a,t_al)$,
$1$ if $(b,k)=(a,t_al+1)$,
and $0$ otherwise.
Then, the claim follows from the
fact that 
$\cR^\nu_l(w(y))
\equiv \tR^\nu(y)$ mod $I_l$,
where $I_l$ is the ideal in Section
\ref{subsec:special},
and that
\begin{align}
\lim_{l\to \infty}\tR^\nu(y)
=
\lim_{l\to \infty}(\tQ\sa_1(y))^{-1}
\frac{\tQ\sa_{t_al+1}(y)}
{\tQ\sa_{t_al}(y)}
=(\tQ\sa_1(y))^{-1}
\end{align}
by the convergence property (\~Q-II).
\end{proof}

We further define series of $y$
\begin{align}
\label{eq:Udef1}
U_a(y)
=y_a \prod_{b=1}^n \tQ_1\sab(y)^{-(\alpha_a|\alpha_b)t_b},
\quad a=1,\, \dots,\, n.
\end{align}

\begin{lem}
\label{lem:k03}
The following equality of series of $y$ holds:
\begin{align}
\label{eq:k07}
\tK^0(y)
=\det_{1\leq a,b \leq n}
\left(
\frac{y_b}{U_a}
\frac{\partial U_a}{\partial y_b}(y)
\right)
\prod_{a=1}^n
\tQ_1\sa(y).
\end{align}
\end{lem}

\begin{proof}
We define series of $y$
\begin{align}
\label{eq:udef1}
u_a(y)=y_a \prod_{b=1}^n(f_b(y))^{(\alpha_a|\alpha_b)t_b},
\quad
a=1,\, \cdots,\, n,
\end{align}
which depend also on $l$.
By Lemma \ref{lem:flim1},
(\ref{eq:k07}) is the limit $l\rightarrow \infty$
of the formula
\begin{align}
\label{eq:k03}
\cK^0_l(w(y))
=\det_{1\leq a,b\leq n}
\left(
\frac{y_b}{u_a}
\frac{\partial u_a}{\partial y_b}(y)
\right)
\prod_{a=1}^n
(f_a(y))^{-1}.
\end{align}
To prove (\ref{eq:k03}),
we use
\begin{align}
\label{eq:dx2}
y_a\frac{\partial}{\partial y_a}
&=
\sum_{m=1}^{t_al} m w\am
\frac{\partial}{\partial w\am},\\
\label{eq:detf3}
\det_{H_l}(\delta_{ab}\delta_{mk}+m\alpha_{abk})
&=
\det_{1\leq a,b\leq n}
(\delta_{ab}+\sum_{k=1}^{t_bl}k\alpha_{abk}),
\quad\text{$\alpha_{abk}$: arbitrary},
\end{align}
where (\ref{eq:detf3})
 is easily shown by elementary transformations.
Let
\begin{align}
F_a(y)=\prod_{b=1}^n (f_b(y))^{(\alpha_a|\alpha_b)t_b}.
\end{align}
By (\ref{eq:K3}),
(\ref{eq:k03}) is equivalent to
the following equality:
\begin{alignat*}{2}
\det_{H_l}\Bigl(
\frac{w\bk}{z\am}
\frac{\partial z\am}{\partial w\bk}
\Bigr)
&=
\det_{H_l}\Bigl(
\delta_{ab}\delta_{mk}+mw\bk
\frac{\partial }{\partial w\bk}
\log F_a
 \Bigr)
&\quad& \text{(by  (\ref{eq:wvrel2}))}\\
&=
\det_{1\leq a,b\leq n}
\Bigl(
\delta_{ab}
+\sum_{k=1}^{t_bl} k w\bk
\frac{\partial }{\partial w\bk}
\log F_a
\Bigr)
&& \text{(by  (\ref{eq:detf3}))}\\
&=
\det_{1\leq a,b\leq n}
\Bigl(
\delta_{ab}
+y_b\frac{\partial }{\partial y_b}
\log F_a
\Bigr)
&&
 \text{(by (\ref{eq:dx2}))}
\\
&=
\det_{1\leq a,b\leq n}
\Bigl(
\frac{y_b}{u_a}
\frac{\partial u_a}{\partial y_b}
\Bigr)
&& \text{(by  (\ref{eq:udef1})).}
\end{alignat*}
\end{proof}

\begin{prop}\label{prop:k01}
The following equality of Laurent
series of $x$ holds:
\begin{align}\label{eq:k01}
K^0(x)=
\det_{1\leq a,b\leq n}
\left(\frac{\partial{Q_1\sa}}{\partial{x_b}}(x)\right).
\end{align}
\end{prop}

\begin{proof}
We recall that
$y_a(x)=\prod_{b=1}^n x_b^{-(\alpha_a|\alpha_b)t_b}$,
$K^0(x)=\tK^0(y(x))$, and 
$Q\sa_1(x)=x_a \tQ\sa_1(y(x))$.
Then, by (\ref{eq:Udef1}), we have
\begin{align}
U_a(y(x))=\prod_{b=1}^n Q\sab_1(x)^{-(\alpha_a|\alpha_b)t_b}.
\end{align}
Therefore, by Lemma \ref{lem:k03},
\begin{align}
K^0(x)=
\det_{1\leq a,b\leq n}
\left(\frac{x_b}{Q\sa_1}
\frac{\partial Q\sa_1}{\partial x_b}(x)
\right)
\prod_{a=1}^n \tQ\sa_1(y(x))
=
\det_{1\leq a,b\leq n}
\left( 
\frac{\partial Q\sa_1}{\partial x_b}(x)
\right).
\end{align}
\end{proof}
\par
{\em Step 2.}
So far, we have not used the assumption of Proposition
\ref{prop:weyl1} that $Q\am$ are $\cW$-invariant yet.

\begin{prop}\label{prop:qx5}
If $Q_1^{(1)}(x)$, \dots, $Q_1^{(n)}(x)$ are $\cW$-invariant,
then
\begin{equation}
\det_{1\leq a,b\leq n}
\left(\frac{\partial Q_1\sa}{\partial x_b}(x)\right)=
\prod_{\alpha \in \Delta_+} (1-e^{-\alpha}).
\end{equation}
\end{prop}
\begin{proof}
It is well-known that
\begin{align*}
e^\rho \prod_{\alpha \in \Delta_+} (1-e^{-\alpha})
= \sum_{w\in \cW} \mathrm{sgn}(w) e^{w(\rho)},
\end{align*}
and the RHS is characterized by (i) the coefficient
of $e^\rho$ is 1; (ii) it is skew $\cW$-invariant.
Therefore, it is enough to show that 
$e^\rho (\partial Q_1/\partial x)(x)$ satisfies the
same properties.
The property (i) follows from the fact $\tQ_1\sa(y)
=1+O(y)$. 
Under the assumption that $Q_1\sa(x)$ are $\cW$-invariant,
the property (ii) is equivalent to the fact
that $e^{-\rho} dx_1\wedge\cdots\wedge dx_n$ is
skew $\cW$-invariant,
This is easily seen by using the following well-known 
transformation property under the simple reflection $s_a$:
\begin{align*}
s_a(e^{-\rho})=(y_a)^{-1}e^{-\rho},
\qquad
s_a(dx_1\wedge \cdots \wedge dx_n)=-y_a
(dx_1\wedge \cdots \wedge dx_n).
\end{align*}
\end{proof}
Proposition \ref{prop:weyl1} follows
from Propositions \ref{prop:k01} and \ref{prop:qx5}.


\end{document}